\documentclass[preprint]{elsarticle}

\usepackage[utf8x]{inputenc}
\usepackage[T1]{fontenc}
\usepackage{amsmath,amssymb,amsthm}
\usepackage[left=34mm,right=34mm]{geometry}
\usepackage{fixmath}
\usepackage{bm}
\usepackage{bbm} %for \mathbbm{1}
\usepackage{graphicx}
\usepackage{xcolor}
\usepackage{booktabs}
\usepackage{listings}

\let\counterwithin\relax
\usepackage{chngcntr}
\usepackage{lmodern}
\usepackage{multirow}
\usepackage{subfig}
\usepackage{placeins}
\usepackage{multirow}
\PassOptionsToPackage{hyphens}{url}
\usepackage{hyperref}
\usepackage[colorinlistoftodos,prependcaption,textsize=tiny]{todonotes}

\graphicspath{{fig/}}

\definecolor{codegreen}{rgb}{0,0.6,0}
\definecolor{codegray}{rgb}{0.5,0.5,0.5}
\definecolor{codepurple}{rgb}{0.58,0,0.82}
\definecolor{backcolour}{rgb}{1,1,1}
\definecolor{codeblue}{RGB}{0,160,238}
\definecolor{codeorange}{RGB}{217,83,25}
\lstdefinestyle{mystyle}{
    language=C++,
    backgroundcolor=\color{backcolour},
    commentstyle=\color{codeblue},
    keywordstyle=\color{codeorange},
    numberstyle=\tiny\color{codegray},
    breakatwhitespace=false,
    breaklines=true,
    captionpos=b,
    keepspaces=true,
    numbers=left,
    numbersep=5pt,
    showspaces=false,
    showstringspaces=false,
    showtabs=false,
    tabsize=2,
    basicstyle={\scriptsize\ttfamily},
    xleftmargin=.05\textwidth,
    xrightmargin=.05\textwidth,
    morekeywords={omp,offload,target,parallel,simd,critical,atomic,linear,declare,uniform,reduction,offload_transfer},
    framextopmargin=2pt,
    frame=tb,
    escapeinside={(|}{|)},
    moredelim=[is][\color{codeorange}]{|*}{*|},
    showlines=true
}
\newcommand{\cpp}{C\nolinebreak\hspace{-.05em}\raisebox{.4ex}{\tiny\bf +}\nolinebreak\hspace{-.10em}\raisebox{.4ex}{\tiny\bf +}}

\newcommand{\vect}[1]{\ensuremath{{\bm{#1}}}}

\newcommand{\va}{{\vect{a}}}
\newcommand{\vb}{{\vect{b}}}
\newcommand{\vc}{{\vect{c}}}
\newcommand{\vd}{{\vect{d}}}

\newcommand{\vf}{{\vect{f}}}
\newcommand{\vg}{{\vect{g}}}

\newcommand{\vn}{{\vect{n}}}

\newcommand{\vs}{{\vect{s}}}

\newcommand{\vw}{{\vect{w}}}
\newcommand{\vx}{{\vect{x}}}
\newcommand{\vy}{{\vect{y}}}

\newcommand{\vkappa}{{\vect{\kappa}}}
\newcommand{\vlambda}{{\vect{\lambda}}}
\newcommand{\vmu}{{\vect{\mu}}}
\newcommand{\vnu}{{\vect{\nu}}}
\newcommand{\vxi}{{\vect{\xi}}}
\newcommand{\vzeta}{{\vect{\zeta}}}

\newcommand{\vzero}{{\vect{0}}}
\newcommand{\vone}{{\vect{1}}}

\newcommand{\mat}[1]{\ensuremath{{\mathsf{#1}}}}

\newcommand{\matK}{{\mat{K}}}

\newcommand{\matM}{{\mat{M}}}

\newcommand{\matV}{{\mat{V}}}

\newcommand{\dif}{\ensuremath{{\mathrm{d}}}}

\DeclareMathOperator{\dist}{dist}

\DeclareMathOperator{\child}{child}
\DeclareMathOperator{\parent}{par}

\DeclareMathOperator{\geomcenter}{center}
\DeclareMathOperator{\level}{level}

\numberwithin{equation}{section}
\numberwithin{figure}{section}
\numberwithin{table}{section}
\AtBeginDocument{\counterwithin{lstlisting}{section}}

\theoremstyle{remark}

\usepackage{algorithm}
\usepackage[noend]{algpseudocode}
\algrenewcommand{\algorithmiccomment}[1]{$\triangleright$ \emph{#1}}

\algblockdefx[Pragma]{Pragma}{EndPragma}[1][]{Start #1}{}
\algnotext[Pragma]{EndPragma}

\usepackage{tikz}
\usepackage{pgfplots}
\pgfplotsset{compat=1.15}

\usetikzlibrary{shadows}
\usetikzlibrary{patterns}
\usetikzlibrary{trees}
\usetikzlibrary{decorations.pathreplacing} %for curly braces
\usepackage{xifthen} %for ifthenelse in tikz

\usepackage{subfig} %for subfloat in figures

\journal{a journal}

\begin{document}

\begin{frontmatter}

\title{A parallel fast multipole method for a space-time boundary element method for the heat equation}

\author[tugraz]{Raphael Watschinger}

\author[it4iaddress,dept470]{Michal Merta\texorpdfstring{\corref{mycorrespondingauthor}}{}}
\cortext[mycorrespondingauthor]{Corresponding author}
\ead{michal.merta@vsb.cz}

\author[tugraz]{Günther Of}

\author[it4iaddress,dept470]{Jan Zapletal}

\address[tugraz]{Institute of Applied Mathematics, Graz University of Technology. \\ Steyrergasse 30, A-8010 Graz, Austria}

\address[it4iaddress]{IT4Innovations, VSB -- Technical University of Ostrava. \\ 17.~listopadu 2172/15, 708~00 Ostrava-Poruba, Czech Republic}

\address[dept470]{Department of Applied Mathematics, VSB -- Technical University of Ostrava. \\ 17.~listopadu 2172/15, 708~00 Ostrava-Poruba, Czech Republic}

\begin{abstract}
We present a novel approach to the parallelization of the parabolic fast multipole method for a space-time boundary element method for the heat equation. We exploit the special temporal structure of the involved operators to provide an efficient distributed parallelization with respect to time and with a one-directional communication pattern. On top, we apply a task-based shared memory parallelization and SIMD vectorization. In the numerical tests we observe high efficiencies of our parallelization approach.
\end{abstract}


\end{frontmatter}

\section{Introduction}
\label{sec:01_intro}
Space-time methods have become a popular subject of research in recent years, e.g., see the proceedings \cite{LangerSteinbach2019}. Current advances in high performance computing (HPC) have facilitated this trend and both fields of research have benefited mutually.
On the one hand, increasing computing power enables the solution of real-world problems by space-time methods with their huge system matrices. On the other hand, the huge global matrices of space-time methods allow to develop solvers with better parallel scalability and to use the full power of modern HPC resources. Thus the absolute computational times can be reduced. Classical approaches to a parallel solution of time-dependent partial differential equations use some decomposition in space together with time-stepping algorithms which are sequential with respect to time and thus limit the parallelization to the spatial components. Therefore, parallel-in-time algorithms, such as parareal \cite{Lions2001} or space-time parallel multigrid \cite{Gander2016} have gained popularity recently, as they enable parallelization in both spatial and temporal dimensions and restore scalability on large numbers of~CPUs.
Another motivation to use space-time methods is space-time adaptivity. In fact, stable space-time formulations enable local space-time refinement. This can result in significantly smaller systems and can help to reduce computational times. In contrast, classical time-stepping schemes are associated with space-time tensor product meshes which lack this flexibility.
As classical schemes can exploit the uniformity of the system heavily, it is very important to develop efficient adaptive solvers to become competetive with respect to computational times.

Space-time boundary element methods for the solution of the transient heat equation have been known for a relatively long time \cite{ArnoldNoon:1987,costabel:1990,HsiaoSaranen93}. % \cite{Chapko199747,LubichSchneider92,ArnoldNoon89,Tausch09,CostabelSaranen03}
In general, boundary element methods (BEM) describe solutions of partial differential equations by boundary data only. They have certain advantages over volume-based methods (such as the finite element method). Let us mention the reduced dimension, a possibly higher order of convergence, a simpler meshing of complicated geometries, simpler mesh manipulation in the case of optimization problems and the easier handling of moving meshes. They are also beneficial for problems on unbounded domains. In addition, BEM is well suited for parallelization. Because of its high computational intensity and the dense structure of its system matrices, the use of BEM can help to leverage the full potential of modern many-core CPUs equipped with wide SIMD (Single Instruction Multiple Data) registers or GPU accelerators.
On the other hand BEM is harder to implement and its efficient application is typically limited to linear partial differential equations with constant coefficients and linear boundary conditions.

BEM system matrices are in general dense which leads to high computational and memory complexities. Therefore, several fast and data-sparse algorithms have been developed to provide efficient solvers with almost linear complexity.
In the case of the heat equation, e.g., there are algorithms based on Fourier series and FFT \cite{GreengardLin2000,Greengard1990949}, the parabolic FMM \cite{Tausch2007,Tausch2012}, or a fast sparse grid method~\cite{HarbrechtTausch2018}.
The parabolic fast multipole method (pFMM) has originally been described for Nyström discretizations~\cite{Tausch2007,Tausch2012} and has been extended to Galerkin discretizations \cite{MessnerSchanzTausch14,MessnerSchanzTausch15} later on. It is based on a clustering of the computational domain in both space and time and an approximation of interactions in well-separated clusters by truncated series expansions. The resulting method can be seen as a combination of a one-dimensional FMM in time and fast Gauss transforms in space.
Commonly, the related space-time system is solved by some sort of block forward elimination successively, i.e.,~like a sequential time-stepping scheme. While this is beneficial in terms of memory requirements and especially in the case of small numbers of time-steps, the parallelization is limited to the spatial components.

Many publications have been devoted to the efficient implementation and parallelization of the FMM (see, e.g., \cite{Abduljabbar2017,abduljabbar2014asynchronous,Aguetal2014,Cruz2011403,Lashuk2012,YokBar2012}) mainly for particle simulations, but also for the solution of classical spatial boundary integral equations \cite{Abduljabbar2019C245,WanCooBetBar2021}, and, less frequently, space-time boundary integral equations \cite{Takahashi2014}. A parallelization of a standard (i.e.~not fast) Galerkin space-time BEM for the heat equation in two spatial dimensions was considered by the authors in \cite{DohZapOfMer2019}.

We aim at developing adaptive space-time boundary element methods to utilize the aforementioned advantages.
There is a first publication \cite{Gantner22} on BEM adaptivity for the heat equation. The numerical examples show a superior order of convergence of space-time adaptivity.
As an important step towards efficient adaptive space-time BEM solvers
we here present a novel parallel version of a space-time FMM for a Galerkin space-time BEM for the heat equation in three spatial dimensions and tensor-product meshes based on the pFMM. We decided for the pFMM since it is well established and seems to be advantageous with respect to parallelization and adaptivity. We enable parallelization in time by solving the whole space-time system at once. Instead of simply transferring one of the parallelization approaches of spatial FMMs to our setting, we use a special task-based concept to tailor our parallel algorithm to the specific structure of the~pFMM. The presented method employs a distribution of the space-time cluster tree with respect to time among MPI (Message Passing Interface) processes. In particular, we exploit the causality of the operators which leads to a one-directional communication in the temporal component. In addition, the computation on individual processes is parallelized in shared memory using OpenMP tasks with explicitly stated dependencies \cite{OpenMP}.

Our parallelization approach in shared and distributed memory is based on a data driven model, instead of a bulk-synchronous parallelization often used in scientific codes.
We use a two-level task-based concept for the parallelization of the space-time FMM.
We group FMM operations as tasks along the underlying temporal tree which reflects the communication pattern of our distribution. This provides a first coarse granularity to represent the dependencies of the FMM data/operations and to steer the load balancing and the MPI parallelization.
By executing tasks (grouped FMM operations) based on individual dependencies, we overcome the strict separation of phases of classical FMM and can arrange the grouped operations more flexibly. This allows us to hide communication and to fill eventual idle times by independent tasks.
When it comes to the execution of temporal tasks on a single MPI process we generate a fine granularity by creating multiple tasks for the related space-time operations by OpenMP. Here the OpenMP tasks scheduler can serve as a buffer for the created OpenMP tasks keeping the cores busy during the computation.

Our numerical experiments cover SIMD vectorization, shared memory performance using up to 36 cores, and scalability tests of the distributed memory parallelization with up to 256 MPI processes (6144 cores). In these tests we observe high efficiencies of our parallelization approach.
We think that such a good scalability cannot be reached for the same examples by some purely spatial parallelization.

The structure of this paper is as follows. In Section~\ref{sec:02_bie_bem} we briefly describe the boundary integral formulation of the Dirichlet initial boundary value problem of the heat equation and its discretization using space-time BEM. We provide a comprehensive description of the space-time FMM for a Galerkin BEM in Section~\ref{sec:03_stFMM}. This serves as a basis for the detailed presentation of our parallelization concept in Section~\ref{sec:04_parallelization}. Finally, numerical experiments and conclusions are provided in Section~\ref{sec:05_experiments} and Section~\ref{sec:06_conclusion}.

\section{Boundary integral formulation of the heat equation and its discretization}
\label{sec:02_bie_bem}
We consider the Dirichlet boundary value problem of the transient heat equation with the heat capacity constant $\alpha>0$ and zero initial condition for a bounded Lipschitz domain~$\Omega \subset \mathbb{R}^3$ with a boundary $\Gamma= \partial \Omega$ as a model problem:
\begin{alignat*}{2}
  \frac{\partial}{\partial t} u(\vx,t) -  \alpha \Delta u(\vx,t) &= 0
  &\quad &\text{for } (\vx,t) \in \Omega \times (0,T),\\
  u(\vx,0) &= 0 &\quad &\text{for } \vx \in \Omega,\\
  u(\vx,t) &= g(\vx,t) &\quad &\text{for } (\vx,t)\in \Sigma := \Gamma \times (0,T).
\end{alignat*}
The solution of this problem can be described by the representation formula
\begin{equation*}
  u(\vx,t) = \widetilde Vw(\vx,t) - Wu(\vx,t) \qquad \text{for } (\vx,t) \in \Omega \times (0,T)
\end{equation*}
with the conormal derivative $w := \alpha\frac{\partial u}{\partial \vn}$, the single layer potential
\begin{equation*}
  \widetilde V w (\vx,t) := \int_0^t \int_{\Gamma} G_\alpha(\vx-\vy,t-\tau) w(\vy,\tau) \,\dif\vs_\vy \,\dif\tau,
\end{equation*}
the double layer potential
\begin{equation*}
  W u (\vx,t) := \int_0^t \int_{\Gamma} \alpha\frac{\partial G_\alpha}{\partial \vn_\vy}(\vx-\vy,t-\tau) u(\vy,\tau) \,\dif\vs_\vy \,\dif\tau,
\end{equation*}
and the heat kernel
\begin{equation} \label{eq:heat_kernel}
  G_\alpha(\vx-\vy,t-\tau) =
  \begin{cases}\displaystyle
   (4\pi\alpha(t-\tau))^{-3/2}\exp\left(-\frac{|\vx-\vy|^2}{4\alpha(t-\tau)}\right)
    & \text{for } \tau<t,\\
    0  & \text{for } \tau>t.
  \end{cases}
\end{equation}
While the Dirichlet datum $u=g$ is given on $\Sigma$, the unknown Neumann datum $w$ can be determined from the boundary integral equation
\begin{equation}
  Vw(\vx,t) = \bigg(\frac{1}{2} I + K \bigg) g(\vx,t) \qquad\text{for almost all } (\vx,t) \in \Sigma
\end{equation}
with the single and double layer boundary integral operators
\begin{align*}
  Vw(\vx,t) &= \int_0^t \int_{\Gamma} G_\alpha(\vx-\vy,t-\tau) w(\vy,\tau) \,\dif\vs_\vy \,\dif\tau, \\
  Ku(\vx,t) &= \int_0^t \int_{\Gamma} \alpha\frac{\partial G_\alpha}{\partial \vn_\vy}(\vx-\vy,t-\tau) u(\vy,\tau) \,\dif\vs_\vy \,\dif\tau.
\end{align*}

As the Neumann datum $w$ is unknown, we compute some numerical approximation. Typically a tensor product mesh is considered for uniform time-steps, $t_{j_t}= j_t h_t$, $j_t=0,\ldots,E_t$, and a fixed spatial surface mesh $\{\gamma_{j_\vx}\}_{j_\vx=1}^{E_\vx}$ of triangles. For such a decomposition $\Sigma_h$ of the lateral boundary ${\Sigma = \Gamma \times (0,T)}$ into space-time boundary elements $\sigma_{j_t,j_\vx} = \gamma_{j_\vx} \times (t_{j_t-1},t_{j_t})$,
a standard approximation of $w$ is given by a piecewise constant approximation
\begin{equation*}
 w_h(\vy,\tau) = \sum_{j_t=1}^{E_t} \sum_{j_\vx=1}^{E_\vx} w_{j_t,j_\vx} \varphi_{j_t,j_\vx}^{0,0}(\vy, \tau)
\end{equation*}
with basis functions $\varphi_{j_t,j_\vx}^{0,0}$ which are one on a space-time boundary element $\sigma_{j_t,j_\vx}$  and zero otherwise.

To find the yet unknown coefficients $w_{j_t,j_\vx}$ we consider the Galerkin variational formulation
\begin{equation*}
\int_{t_{k_t-1}}^{t_{k_t}} \int_{\gamma_{k_\vx}} Vw_h(\vx,t) ds_\vx dt
=  \int_{t_{k_t-1}}^{t_{k_t}} \int_{\gamma_{k_\vx}} \bigg(\frac{1}{2} I + K \bigg) g_h(\vx,t) \dif\vs_\vx \,\dif t
\end{equation*}
for all $k_t=1,\ldots,E_t$ and $k_\vx=1,\ldots,E_\vx$. Note that we have replaced the Dirichlet datum~$g$ by a $L_2(\Sigma)$ projection~$g_h$, which is piecewise constant in time, but piecewise linear and globally continuous in space.
The equivalent system of linear equations to find $w_h$ or rather the vector $\vw$ of coefficients~$w_{j_t,j_\vx}$ is
\begin{equation}\label{02:linear:system}
  \matV_h \vw = \bigg(\frac{1}{2} \matM_h + \matK_h \bigg) \vg,
\end{equation}
where the matrix entries are given by
\begin{equation*}
  \mat{V}[(k_t-1)E_\vx +k_\vx,(j_t-1)E_\vx +j_\vx] =
\int_{t_{k_t-1}}^{t_{k_t}} \int_{\gamma_{k_\vx}}\int_{t_{j_t-1}}^{t_{j_t}}\int_{\gamma_{j_\vx}} G_\alpha(\vx-\vy,t-\tau) \,\dif\vs_\vy \,\dif\tau \,\dif\vs_\vx \,\dif t
\end{equation*}
for $k_\vx,j_\vx=1,\ldots,E_\vx$ and $k_t,j_t=1,\ldots,E_t$. Due to the causality of the heat kernel this matrix has lower triangular block structure.
$\mat{K}_h$ is defined similarly but for trial functions which are piecewise constant in time, piecewise linear and globally continuous in space. $\mat{M}_h$ denotes the related mass matrix.
Please check \cite{zapletal2021semianalytic} for details on the discretization and implementation. Some detailed analysis of the integral equations and the presented boundary element method is provided in \cite{costabel:1990,DohNiiSte2019}.

Note that system \eqref{02:linear:system} is huge and that the matrices $\mat{V}_h$ and $\mat{K}_h$ are dense except for their lower triangular block structure. Thus a standard BEM is limited to small problems. In general a data-sparse method such as the FMM is necessary to solve large-scale problems.

\section{The sequential space-time FMM algorithm for the heat equation}
\label{sec:03_stFMM}
Following \cite{Tausch2007,Tausch2012, MessnerSchanzTausch14} we first present a sequential space-time FMM for the multiplication of a vector $\vw$ by the single layer operator matrix $\mat{V}_h$ before dealing with its parallelization. The considered FMM is based on a suitable expansion of the heat kernel~$G_\alpha$ and a clustering of the computational domain, and can also be applied to other BEM matrices of the heat equation like~$\mat{K}_h$ in \eqref{02:linear:system} with slight modifications. In the parabolic fast multipole method in \cite{Tausch2007,Tausch2012, MessnerSchanzTausch14} the matrix-vector multiplication and solution of linear systems is executed in a forward-sweeping manner, based on the causality of the operators. However, this does not allow for a parallelization in time, which is our aim. Thus, we consider the method more like a standard FMM realizing the full matrix-vector product at once, and use an iterative solver like GMRES for system~\eqref{02:linear:system}. We still want to highlight the special temporal structure of the space-time method in our description, since it forms the basis of our parallelization strategy.

\subsection{A separable approximation of the heat kernel} \label{subsec:kernel_approx}

\begin{figure}
  \centering
  \subfloat[]{\includegraphics{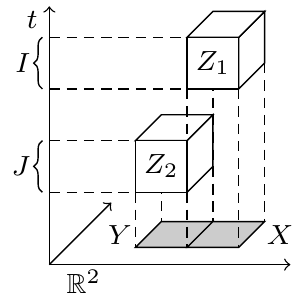}\label{fig:boxes:a}}
  \qquad \qquad \qquad
  \subfloat[]{\includegraphics{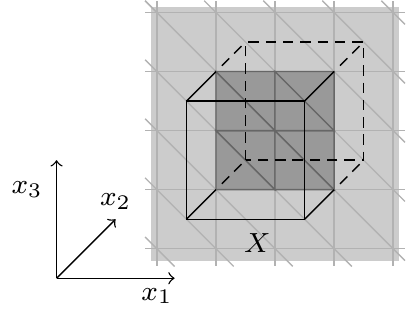}\label{fig:boxes:b}}
  \caption{(a) shows a space-time target box $Z_1 = X \times I$ and a source box $Z_2 = Y \times J$. For the sake of simplicity the boxes are drawn in 3D instead of 4D. In (b) a box $X$ in three spatial dimensions is shown that represents the spatial part of a space-time box. It contains some triangles (dark gray) of a spatial surface mesh $\{\gamma_{j_\vx}\}_{j_\vx}^{E_\vx}$ (all gray triangles) as it is the case for the boxes constructed in Section \ref{sec:4d:tree}.} \label{fig:boxes}
\end{figure}

Throughout the section we regard the heat kernel as a function of the differences $\vx-\vy$ and~$t-\tau > 0$ as in \eqref{eq:heat_kernel}. We restrict the variables $(\vx,t)$ to a 4D target box $Z_1 = X \times I$ and $(\vy,\tau)$ to a 4D source box $Z_2 = Y \times J$, where $I=(c_4,c_4+2\tilde{h}_t]$ and $J=(d_4,d_4+2\tilde{h}_t]$ are two intervals of length $2\tilde{h}_t$ such that $c_4 > d_4$ and $\dist(I,J) = c_4 - d_4 - 2\tilde{h}_t > 0$, and $X = (\vc, \vc + 2\tilde{h}_x \vone]$ and $Y = (\vd, \vd + 2\tilde{h}_x \vone]$ are cubes in $\mathbb{R}^3$ with edge length~$2 \tilde{h}_x$. Here we use the notation
\begin{equation*}
  (\va, \vb] := (a_1,b_1] \times (a_2,b_2] \times (a_3,b_3].
\end{equation*}
A sketch of the boxes $Z_1$ and $Z_2$ is given in Figure \ref{fig:boxes:a}. For~$(\vx,t) \in Z_1$ and $(\vy,\tau) \in Z_2$ there holds $t - \tau \geq \dist(I,J) > 0$ and thus the heat kernel is smooth. As in \cite{Tausch2007, Tausch2012, MessnerSchanzTausch14} we interpolate it in the temporal intervals $I$ and $J$ and approximate it in the spatial boxes $X$ and $Y$ by means of a truncated Chebyshev expansion.

For this purpose, let $T_k(x) = \cos(k \arccos(x))$ be the Chebyshev polynomials of order $k$ on~$[-1,1]$ and let $\{\xi_k^{(m_t)}\}_{k=0}^{m_t}$ be the Chebyshev nodes of order $m_t+1$ on $[-1,1]$, i.e.,~the roots of~$T_{m_t+1}$. On an interval $I = (a,b]$ we consider the transformed points ${\xi_{I,k}^{(m_t)} := \varphi_{I}(\xi_k^{(m_t)})}$, where~$\varphi_I$ is the affine transformation from $(-1,1]$ to $I$, and the associated Lagrange polynomials
\begin{equation*}
  L_{I,b}(t) := \prod_{k \neq b}\frac{t - \xi_{I,k}^{(m_t)}}{\xi_{I,b}^{(m_t)}-\xi_{I,k}^{(m_t)}}
    \qquad \text{for all } b \in\{1, ..., m_t+1\}.
\end{equation*}
For the expansion in space we use the transformed Chebyshev polynomials $T_{(a,b],k}:=T_k \circ \varphi_{(a,b]}^{-1}$ on intervals $(a,b]$ and their tensor products $T_{(\va,\vb],\vkappa}(\vx) := \prod_{j} T_{(a_j,b_j],\kappa_j}(x_j)$, where $\vkappa$ is a multi-index in~$\mathbb{N}_0^3$. By interpolating $G_\alpha$ in the temporal points $\{\xi_{I,k}^{(m_t)}\}_k$ and $\{\xi_{J,k}^{(m_t)}\}_k$ and approximating the resulting function via a truncated Chebyshev expansion in $X$ and $Y$ we~get
\begin{equation} \label{eq:kernel_approx}
  G_\alpha(\vx-\vy,t-\tau)
  \approx \sum_{a,b=0}^{m_t} \sum_{|\vkappa + \vnu| \leq m_\vx}
    E_{a,\vkappa,b,\vnu} T_{X,\vnu}(\vx) T_{Y,\vkappa}(\vy) L_{I,b}(t) L_{J,a}(\tau)
\end{equation}
in $Z_1 \times Z_2$. Here, $m_\vx \geq 0$ is the expansion order in space, and
\begin{equation} \label{eq:exp_coeffs_3d}
  E_{a,\vkappa,b,\vnu} = \frac{1}{(4\pi\alpha (\xi_{I,b}^{(m_t)}-\xi_{J,a}^{(m_t)}))^{3/2}}
    \prod_{j=1}^3 E_{\kappa_j,\nu_j}(r_j,d_{a,b})
\end{equation}
are the expansion coefficients, where
\begin{align}
  r_j &:= (c_j -d_j)/\tilde{h}_x \label{eq:def_r_j}, \\
  d_{a,b} &:= 4 \alpha (\xi_{I,b}^{(m_t)}-\xi_{J,a}^{(m_t)}) / \tilde{h}_x^2, \label{eq:def_d_ab} \\
  E_{k,\ell}(r,d_{a,b})
  &= \frac{\lambda_k \lambda_l}{({m_\vx} + 1)^2} \sum_{n,m=0}^{m_\vx}
    \exp\left( -\frac{|r+\xi_n^{({m_\vx})}-\xi_m^{({m_\vx})}|^2}{ d_{a,b}} \right)
    T_\ell(\xi_n^{({m_\vx})})T_k(\xi_m^{({m_\vx})}) \label{eq:exp_coeffs_1d}
\end{align}
and $\lambda_0 = 1$, $\lambda_k = 2$ for all $k>0$, cf.~\cite[Section~5.3, page 209]{Tausch2012}.

Let us comment on the approximation quality of \eqref{eq:kernel_approx}. The temporal interpolation error converges exponentially to zero with respect to the interpolation degree $m_t$ if the time intervals are well-separated, see e.g.,~\cite[Lemma 4.1 and Equation~(4.42)]{Mes2014}. For the truncated Chebyshev expansion one can show super-exponential convergence of the approximation error without requiring a separation of the spatial boxes $X$ and~$Y,$ see \cite[Section~5.3, page 209]{Tausch2012}. However, the effective approximation quality suffers for small values of $d_{a,b}$. Therefore, we bound $d_{a,b}$ from below in the later application by choosing the spatial box half-size~$\tilde{h}_x$ for a given temporal interval half-size~$\tilde{h}_t$ such that
\begin{equation} \label{eq:rel_h_t_h_x}
    \frac{\tilde{h}_x^2}{4 \alpha\, \tilde{h}_t} \leq c_{\mathrm{st}}
\end{equation}
for some constant $c_{\mathrm{st}} > 0$ \cite[cf. $\rho$ in (27) and Section 5.4]{Tausch2012}. Since $(\xi_{I,b}^{(m_t)}-\xi_{J,a}^{(m_t)}) \geq \dist(I,J)$, \eqref{eq:rel_h_t_h_x} implies that $d_{a,b} \gtrsim c_{\mathrm{st}}^{-1}$ if we guarantee that $\dist(I,J) \gtrsim \tilde{h}_t$, i.e.,~$I$ and $J$ are well-separated.

An important observation is that the heat kernel decays exponentially in space for fixed temporal variables. Therefore, we do not have to approximate the values of the kernel in boxes $Z_1 = X \times I$ and $Z_2 = Y \times J$ but can instead neglect them if the distance of the spatial boxes $X$ and $Y$ is large compared to the distance of the time intervals $I$ and $J$.

\subsection{A 4D space-time box cluster tree}\label{sec:4d:tree}
\begin{figure}
  \centering
  \subfloat[]{\includegraphics{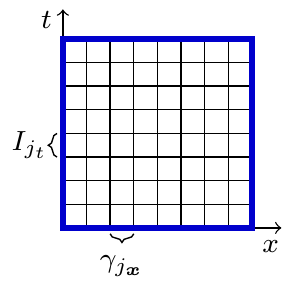}\label{fig:box_cluster_tree:a}}
  \quad
  \subfloat[]{\includegraphics{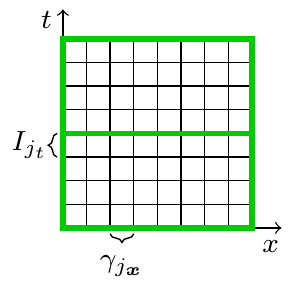}\label{fig:box_cluster_tree:b}}
  \quad
  \subfloat[]{\includegraphics{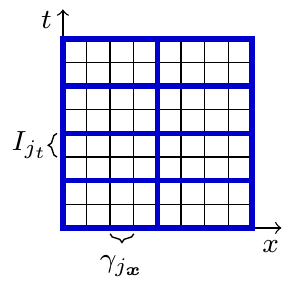}\label{fig:box_cluster_tree:c}}
  \quad
  \subfloat[]{\includegraphics{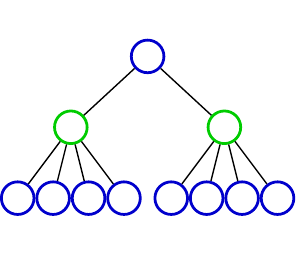}\label{fig:box_cluster_tree:d}}
  \caption{Construction of a space-time cluster tree in 2D. Instead of 4D space-time elements $\sigma$ we consider two-dimensional rectangular elements consisting of temporal intervals $\{J_{j_t}\}_{j_t=1}^{E_{t}}$ and spatial 1D elements $\{\gamma_{j_\vx}\}_{j_\vx=1}^{E_\vx}$. First, a box $Z^{(0)}$ at level zero is constructed that contains all these elements, see~(a). This box is refined recursively as described in Algorithm~\ref{alg:box_cluster_tree}. A purely temporal refinement leads to the two~boxes at level~one in~(b). Refining these boxes in space and time results in eight boxes at level~two, see~(c). By connecting all child boxes with their respective parent box we obtain the space-time tree in (d).} \label{fig:box_cluster_tree}
\end{figure}

For the FMM algorithm we establish a hierarchy of 4D boxes to partition the space-time tensor mesh $\Sigma_h$ appropriately. The resulting structure is denoted as a box cluster tree~$\mathcal{T}_\Sigma$. Our approach is similar to the one in \cite{Tausch2007,MessnerSchanzTausch14}. However, instead of building separate trees in space and time first, we directly establish a 4D~space-time tree, which is a more general approach applicable also to space-time meshes without a strict tensor product structure. For the construction we use a recursive refinement strategy which is described in Algorithm~\ref{alg:box_cluster_tree} and illustrated in Figure~\ref{fig:box_cluster_tree}.

\begin{algorithm}
  \caption{Construction of a 4D space-time box cluster tree~$\mathcal{T}_\Sigma$ \label{alg:box_cluster_tree}}
  \begin{algorithmic}[1]
  \State Let a space-time tensor mesh $\Sigma_h$ inside a 4D box $Z^{(0)} = (\va,\va+2h_x^{(0)}\vone] \times (0,T]$ be given such that $h_x^{(0)}$ and $h_t^{(0)}:=T/2$ satisfy \eqref{eq:rel_h_t_h_x}. \label{algline_bct_assumptions}
  \State Let a bound $n_{\mathrm{max}}$ for the number of elements in a leaf box and $c_{\mathrm{st}} > 0$ for \eqref{eq:rel_h_t_h_x} be given.
  \State Construct an empty tree $\mathcal{T}_\Sigma$ and add $Z^{(0)}$ as its root.
  \State Call \Call{RefineCluster}{$Z^{(0)}$, $\mathcal{T}_\Sigma$}
  \Statex
  \Function{RefineCluster}{$Z$, $\mathcal{T}_\Sigma$}
    \If{$\#\{ \sigma \subset \Sigma_h : \geomcenter(\sigma) \in Z \} \geq n_{\max} $}
      \State Let $\ell = \level(Z)$, $h_t^{(\ell+1)} = 2^{-\ell-1} h_t^{(0)}$ and $\tilde{h}_x$ be the spatial half-size of $Z$. \label{algline:bct_half_sizes}
      \If{$h_t^{(\ell+1)}$ and $\tilde{h}_x$ satisfy \eqref{eq:rel_h_t_h_x}}
        \State Subdivide $Z$ into $n_\mathrm{C}=2$ children $\{Z_k\}_{k=1}^{n_\mathrm{C}}$ by a temporal refinement. \label{algline:bct_temporal_ref}
      \Else
        \State Subdivide $Z$ into $n_\mathrm{C}=16$ children $\{Z_k\}_{k=1}^{n_\mathrm{C}}$ by a space-time refinement.\label{algline:bct_spacetime_ref}
      \EndIf
      \For{$k=1$, \ldots, $n_\mathrm{C}$}
        \If{$\#\{ \sigma \subset \Sigma_h : \geomcenter(\sigma) \in Z_k \} \neq 0$}
          \State Add $Z_k$ to $\mathcal{T}_\Sigma$ as child of $Z$ and call \Call{RefineCluster}{$Z_k$, $\mathcal{T}_\Sigma$}.
        \EndIf
      \EndFor
    \EndIf
  \EndFunction
  \end{algorithmic}
\end{algorithm}

The refinement of a box $Z=(\va,\vb]\times(c,d]$ in lines \ref{algline:bct_temporal_ref} and \ref{algline:bct_spacetime_ref} of Algorithm~\ref{alg:box_cluster_tree} is done as follows. The temporal part $(c,d]$ is split into halves $(c,\widetilde{c}]$ and $(\widetilde{c},d]$ with $\widetilde{c}=(c+d)/2$ in most cases. Only if $\widetilde{c}$ does not coincide with a time-step $t_k$ from the considered mesh, we choose $\widetilde{c}=t_{k^*}$ as splitting point instead, where $t_{k^*}$ is the time-step closest to the center $(c+d)/2$. In case of a purely temporal refinement, we split~$Z$ into the boxes $Z_1=(\va,\vb]\times(c,\widetilde{c}]$ and $Z_2=(\va,\vb] \times (\widetilde{c},d]$. In case of a space-time refinement we additionally split the spatial part $(\va,\vb]$ of~$Z$ uniformly into 8 boxes $(\va,\widetilde{\va}]$, $\dots$, $(\widetilde{\va},\vb]$ where $\widetilde{\va} = 1/2\,(\va + \vb)$, and get 16 space-time boxes~$\{Z_j\}_{j=1}^{16}$ as combinations of these refined spatial boxes and the temporal intervals $(c,\widetilde{c}]$ and~$(\widetilde{c},d]$. Note that due to \eqref{eq:rel_h_t_h_x} we alternate between purely temporal and space-time refinements in the construction of $\mathcal{T}_\Sigma$, at least after some initial temporal refinements.

With the described splitting in time it is guaranteed that the temporal part $(t_{k_{t}-1},t_{k_t})$ of a space-time element $\sigma_{k_t,k_\vx}=\gamma_{k_\vx} \times (t_{k_{t}-1},t_{k_t})$ is always fully contained in the temporal interval $I$ of a box~$Z=X\times I$ if its center is in~$Z$. For such a space-time element $\sigma_{k_t,k_\vx}$ we want to ensure in addition that $\gamma_{k_\vx} \subset X$ as it is shown in Figure~\ref{fig:boxes:b}. If this is not the case, we pad the boxes in $\mathcal{T}_\Sigma$ appropriately in a post-processing step, i.e.,~we extend the spatial size of a box $Z$ if necessary such that for all~$\sigma$ with $\geomcenter(\sigma) \in Z$ there holds $\sigma \subset Z$. To retain the uniformity of the spatial parts of the boxes at a given level of $\mathcal{T}_\Sigma$, we pad all of them by the same amount in all directions. In addition, we pad boxes at level $\ell$ of $\mathcal{T}$ at least by the same amount as boxes at level $\ell+1$. This ensures that the children of a box $Z$ in $\mathcal{T}_\Sigma$ are fully contained in $Z$, which we also need later on.

A few aspects of Algorithm~\ref{alg:box_cluster_tree} require additional attention. The estimate \eqref{eq:rel_h_t_h_x} does not have to be satisfied in general for given $h_x^{(0)}$, $h_t^{(0)}$ and $c_{\mathrm{st}}$ as required in line~\ref{algline_bct_assumptions}. However, it can be established by additional refinements of the initial spatial box. Due to the non-uniform refinement in time and the padding in space, \eqref{eq:rel_h_t_h_x} might also be violated for other boxes in $\mathcal{T}_\Sigma$, but for suitably regular meshes $\Sigma_h$ it will still hold for a slightly larger constant $c_{\mathrm{st}}$. Finally, we want to point out that the refinement process should be stopped earlier for a box $Z$ if it contains a space-time element $\sigma = \gamma \times (t_{j-1},t_j)$ whose temporal or spatial size is considerably larger than the temporal or spatial half-size of $Z$. This is in particular the case if all of the elements in $Z$ share the same temporal component.

We denote the set of boxes/clusters at level $\ell$ of $\mathcal{T}_\Sigma$ by $\mathcal{T}_\Sigma^{(\ell)}$, its leaves by $\mathcal{L}_\Sigma$ and its depth, which is the largest level attained by any of its clusters, by $p(\mathcal{T}_\Sigma)$. For a cluster $Z \in \mathcal{T}_\Sigma$ we denote the set of all its children by $\child(Z)$ and its parent by $\parent(Z)$. By~$\hat{Z}$ we denote the set of all indices $(k_t,k_\vx)$ such that $\sigma_{k_t,k_\vx} = \gamma_{k_\vx} \times (t_{k_t-1},t_{k_t}) \in Z$.

\subsection{Nearfield and interaction lists of boxes in a space-time cluster tree}\label{subsec:nfia}
In Section \ref{subsec:kernel_approx} we have approximated the heat kernel $G_\alpha(\vx-\vy,t-\tau)$ for~$(\vx,t)$ in a target box~$Z_1$ and $(\vy,\tau)$ in a source box $Z_2$. We recall that for this approximation $Z_1$ and $Z_2$ have to be separated in time but not in space, and that the values of $G_\alpha$ are negligibly small if the spatial distance of $Z_1$ and $Z_2$ is large enough. Based on these observations we define the nearfield and interaction lists of target boxes in $\mathcal{T}_\Sigma$ which will determine the operations in the FMM algorithm.

We start by considering the temporal components of boxes in $\mathcal{T}_\Sigma$. By construction of $\mathcal{T}_\Sigma$ there exist at most $2^\ell$ distinct time intervals that are components of boxes in~$\mathcal{T}_\Sigma^{(\ell)}$. These intervals can be organized in a binary tree $\mathcal{T}_I$. In general, $\mathcal{T}_I$ is a full binary tree with depth $p(\mathcal{T}_I)=p(\mathcal{T}_\Sigma)$, but does not have to be a perfect binary tree. The intervals in $\mathcal{T}_I^{(\ell)}$, i.e.,~at level $\ell$ of $\mathcal{T}_I$, are numbered in ascending order from $0$ to $2^{\ell}-1$ skipping the numbers of potentially missing intervals. The leaves of $\mathcal{T}_I$ are denoted by $\mathcal{L}_I$.

\begin{figure}
  \centering
  \subfloat[]{\includegraphics{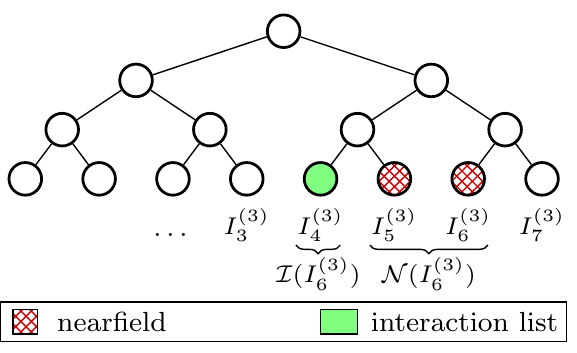}\label{fig:nearfield_and_interaction_lists:a}}
  \qquad
  \subfloat[]{\includegraphics{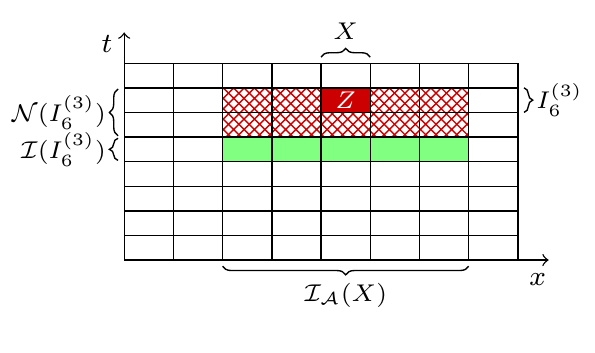}\label{fig:nearfield_and_interaction_lists:b}}
  \caption{Nearfield and interaction lists of temporal intervals and space-time boxes. In (a) a temporal binary tree $\mathcal{T}_I$ is drawn, and the nearfield and interaction list of the interval $I_6^{(3)}$ as defined in \eqref{eq:03:ILtemp} are marked. The figure in (b) shows space-time boxes in 2D. For the box $Z = X \times I_6^{(3)}$ the nearfield and interaction list are highlighted. In addition, the interaction area of $X$ and the nearfield and interaction list of $I_6^{(3)}$ are marked. Note that in this example the interaction area of $X$ contains spatial boxes with grid distance from $X$ bounded by two.} \label{fig:nearfield_and_interaction_lists}
\end{figure}

For a target interval $I_k^{(\ell)}$ in $\mathcal{T}_I^{(\ell)}$ we define the \emph{nearfield} $\mathcal{N}(I_k^{(\ell)})$ and \emph{interaction list} $\mathcal{I}(I_k^{(\ell)})$ by
\begin{align}
  \mathcal{N}(I_k^{(\ell)}) &:=
  \begin{cases}
    \{I_k^{(\ell)}\}, \qquad &\text{if } k = 0, \\
    \left(\{I_{k-1}^{(\ell)}, I_k^{(\ell)}\} \cap \mathcal{T}_I^{(\ell)}\right) \cup \left( \mathcal{N}(\parent(I_k^{(\ell)})) \cap \mathcal{L}_I \right), \qquad &\text{otherwise},
  \end{cases} \nonumber\\
  \mathcal{I}(I_k^{(\ell)}) &:=
  \begin{cases}
    \emptyset, \qquad &\text{if } k \in \{0,1\}, \\
    \{I_{k-2}^{(\ell)}\} \cap \mathcal{T}_I^{(\ell)}, \qquad &\text{if } k \geq 2 \text{ and } k \text{ is even}, \\
    \{I_{k-3}^{(\ell)}, I_{k-2}^{(\ell)}\} \cap \mathcal{T}_I^{(\ell)}, \qquad &\text{if } k \geq 2 \text{ and } k \text{ is odd}.
  \end{cases}\label{eq:03:ILtemp}
\end{align}
Due to the causality of the heat kernel, both sets include only intervals with indices $j\leq k$. Figure~\ref{fig:nearfield_and_interaction_lists:a} gives an example of the nearfield and interaction list of an interval in $\mathcal{T}_I$. We see that the nearfield of an interval $I$ contains directly neighboring intervals including $I$ itself ($\mathcal{N}(I_6^{(3)})=\{I_5^{(3)}, I_6^{(3)}\}$ in Figure~\ref{fig:nearfield_and_interaction_lists:a}). In the general case the nearfield can also include intervals in the nearfield of $I$'s parent which are leaves. The interaction list of $I$ consists of intervals $J$ which are separated from $I$ and whose parents are neighbors of $I$'s parent ($\mathcal{I}(I_6^{(3)}) = \{ I_4^{(3)}\}$ in Figure~\ref{fig:nearfield_and_interaction_lists:a}). Such intervals $I$ and $J$ are thus suitable for the kernel approximation in Section~\ref{subsec:kernel_approx}. Note that earlier time intervals are not contained in the interaction list of a cluster $I$, because they can be handled on a coarser level in the tree in the later FMM algorithm, see Section~\ref{subsec:main_fmm}.

For the spatial component $X^{(\ell)}$ of a box $Z^{(\ell)} \in \mathcal{T}_\Sigma^{(\ell)}$ we want to introduce the \emph{interaction area, which is a certain local neighborhood of $X^{(\ell)}$. Due to the uniform cluster-refinement in space, $X^{(\ell)}$} is contained in a regular grid $\mathcal{G}^{\ell}$ consisting of $8^{\ell_x}$ possibly overlapping boxes, where $\ell_x$ is the number of spatial refinements of boxes in~$\mathcal{T}_\Sigma^{(\ell)}$ and depends on $\ell$. Boxes in the grid $\mathcal{G}^{\ell}$ can be labeled by using multi-indices in $\{0,...,2^{\ell_x}-1\}^3$. We say that two boxes~$X$ and~$Y$ in $\mathcal{G}^{\ell}$ have \emph{grid distance} $n$ if the related multi-indices $\vxi$ and $\vzeta$ satisfy $n = \max_j\{|\xi_j-\zeta_j|\}$. For a fixed parameter $n_{\mathrm{tr}}$ we define the interaction area $\mathcal{I_A}(X^{(\ell)})$ of $X^{(\ell)}$ in~$\mathcal{G}^{(\ell)}$ by
\begin{equation}\label{eq:03:Iarea}
  \mathcal{I_A}(X^{(\ell)}) := \{ Y^{(\ell)} \in \mathcal{G}^{(\ell)} : \text{ the grid distance of } X^{(\ell)} \text{ and } Y^{(\ell)} \text{ is at most } n_{\mathrm{tr}} \}.
\end{equation}

Finally, we define the nearfield and interaction list of a box $Z^{(\ell)} = X \times I \in \mathcal{T}_\Sigma^{(\ell)}$~by
\begin{align}
  \mathcal{N}(Z^{(\ell)}) &:= \left\{ Z^{(\ell)}_{\mathrm{src}} = Y \times J \in \mathcal{T}_\Sigma^{(\ell)} :
    J \in \mathcal{N}(I) \text{ and } Y \in \mathcal{I_A}(X) \right\}
  \cup \left(\mathcal{N}(\parent(Z^{(\ell)})) \cap \mathcal{L}_\Sigma \right), \label{eq:def_st_nearfield} \\
  \mathcal{I}(Z^{(\ell)}) &:= \left\{ Z^{(\ell)}_{\mathrm{src}} = Y \times J \in \mathcal{T}_\Sigma^{(\ell)} : J \in \mathcal{I}(I) \text{ and } Y\in \mathcal{I_A}(X) \right\}. \label{eq:def_st_interaction_list}
\end{align}
The lists in \eqref{eq:def_st_nearfield} and \eqref{eq:def_st_interaction_list} are sketched in Figure~\ref{fig:nearfield_and_interaction_lists:b}. We see that a box $Z_1=Y\times J$ is typically in the nearfield or interaction list of another box $Z_2 = X \times I$ if $Y$ is in the interaction area of $X$ and $J$ is in the nearfield or interaction list of $I$, respectively. Only early leaf clusters in $\mathcal{T}_\Sigma$ have to be treated separately. A box $Z_1=Y\times J$ whose spatial component is not in $\mathcal{I_A}(X)$ is excluded from both lists. This is motivated by the observations in the last paragraph of Section~\ref{subsec:kernel_approx}. Note that the same cutting parameter $n_\mathrm{tr}$ can be chosen for all boxes in $\mathcal{T}_\Sigma$, cf.~\cite[p.210]{Tausch2012}.

\subsection{The main space-time FMM algorithm} \label{subsec:main_fmm}

With the box cluster tree $\mathcal{T}_\Sigma$ and the interaction and nearfield lists of its clusters we construct a partition of the matrix $\mat{V}_h$ into blocks. By $\mat{V}_h|_{\hat{Z}_{\mathrm{tar}} \times \hat{Z}_{\mathrm{src}}}$ we denote the block of $\mat{V}_h$ whose rows correspond to indices $(j_t,j_\vx) \in \hat{Z}_{\mathrm{tar}}$ and columns to indices~$(k_t,k_\vx) \in \hat{Z}_{\mathrm{src}}$. We decompose $\mat{V}_h$ into \emph{admissible blocks} corresponding to indices $\hat{Z}_{\mathrm{tar}} \times \hat{Z}_{\mathrm{src}}$ with $Z_\mathrm{src} \in \mathcal{I}(Z_\mathrm{tar})$, \emph{inadmissible blocks} $\hat{Z}_{\mathrm{tar}} \times \hat{Z}_{\mathrm{src}}$ where $Z_\mathrm{tar} \in \mathcal{L}_\Sigma$ and $Z_\mathrm{src} \in \mathcal{N}(Z_\mathrm{tar})$, and remaining blocks, whose entries are zero due to the lower triangular block structure of $\mat{V}_h$ or negligibly small due to the exponential decay in space of the kernel $G_\alpha$. The FMM algorithm is used to compute an efficient approximation of the product $\mat{V}_h \vw$ using this partition, see Algorithm~\ref{alg:fmm}. In the following we describe the related operations.

Inadmissible blocks of $\mat{V}_h$ are small by construction, since a leaf box $Z_\mathrm{tar} \in \mathcal{L}_\Sigma$ contains only few space-time elements. Hence, we can afford to store and apply them directly, i.e.,~we compute the product ${\widetilde{\vf}|_{\hat{Z}_\mathrm{tar}} = (\mat{V}_h|_{\hat{Z}_\mathrm{tar} \times \hat{Z}_\mathrm{src}} \vw|_{\hat{Z}_\mathrm{src}})}$ as part of $\vf = \mat{V}_h \vw$ by
\begin{equation} \label{eq:eval_sl_nf_block}
  \widetilde{f}_{k_t,k_\vx}
    = \sum_{(j_t,j_\vx) \in \hat{Z}_{\mathrm{src}}} w_{j_t,j_\vx} \int_{t_{k_t-1}}^{t_{k_t}} \int_{\gamma_{k_\vx}}
      \int_{t_{j_t-1}}^{t_{j_t}} \int_{\gamma_{j_\vx}}
      G_\alpha(\vx-\vy,t-\tau)
      \,\dif\vs_\vy \,\dif\tau \,\dif\vs_\vx \,\dif t
\end{equation}
for all $(k_t, k_\vx) \in \hat{Z}_{\mathrm{tar}}$, where $w_{j_t,j_\vx}$ denotes a coefficient of the vector $\vw$. The computation of the corresponding integrals is discussed in \cite{zapletal2021semianalytic}.

For an admissible block $\mat{V}_h|_{\hat{Z}_\mathrm{tar} \times \hat{Z}_\mathrm{src}}$ with boxes $Z_\mathrm{tar} = X \times I$ and $Z_\mathrm{src} = Y \times J$ we can replace the kernel $G_\alpha$ in~\eqref{eq:eval_sl_nf_block} by its approximation in~\eqref{eq:kernel_approx}. The result can be computed in 3 steps:

\noindent \emph{S2M}:
For $a \in \{0, ..., m_t\}$ and $\vkappa \in \mathbb{N}_0^3$ with $|\vkappa|\leq m_\vx$ compute the \emph{moments} $\vmu(Z_\mathrm{src})$~by
\begin{equation} \label{eq:s2m}
  \mu_{a,\vkappa}(Z_\mathrm{src})
  := \sum_{(j_t,j_\vx) \in \hat{Z}_{\mathrm{src}}} w_{j_t,j_\vx}
    \int_{t_{j_t-1}}^{t_{j_t}} \int_{\gamma_{j_\vx}} T_{Y,\vkappa}(\vy) L_{J,a}(\tau) \,\dif\vs_\vy \,\dif\tau.
\end{equation}

\noindent \emph{M2L}:
For $b \in \{0, ..., m_t\}$ and $\vnu \in \mathbb{N}_0^3$ with $|\vnu| \leq m_\vx$ compute the \emph{local contributions}~$\vlambda(Z_\mathrm{tar})$ by
\begin{equation} \label{eq:m2l}
  \lambda_{b, \vnu}(Z_\mathrm{tar})
  := \sum_{a=0}^{m_t} \sum_{|\vkappa + \vnu| \leq m_\vx} E_{a,\vkappa,b,\vnu}\, \mu_{a, \vkappa}(Z_\mathrm{src}).
\end{equation}

\noindent \emph{L2T}:
For all $(k_t, k_\vx)$ in $\hat{Z}_{\mathrm{tar}}$ evaluate
\begin{equation} \label{eq:l2t}
  \widetilde{f}_{k_t,k_\vx} = \sum_{b=0}^{m_t} \sum_{|\vnu|\leq m_\vx}
    \lambda_{b, \vnu}(Z_\mathrm{tar})
    \int_{t_{k_t-1}}^{t_{k_t}} \int_{\gamma_{k_\vx}} T_{X,\vnu}(\vx) L_{I,b}(t) \,\dif\vs_\vx \,\dif t.
\end{equation}

To enhance the performance we additionally use a nested computation of the moments and local contributions, see e.g.,~\cite[Sections~4.2 and~4.3]{Tausch2007}. Moments of a non-leaf cluster are computed from the moments of its children via M2M operations. Local contributions of a non-leaf cluster are passed down to its children with an L2L operation, and evaluated together with the children's local contributions. We distinguish temporal and space-time M2M and L2L operations. Note that the temporal L2L and space-time L2L operations are just the transposed operations of the corresponding M2M operations, which is why we discuss only the latter.

\noindent \emph{Temporal M2M}:
For a box $Z_\mathrm{p} = X \times I_\mathrm{p}$ whose children are refined only in time there holds
\begin{equation} \label{eq:temp_m2m}
  \mu_{a_\mathrm{p}, \vkappa}(Z_\mathrm{p})
  = \sum_{\substack{Z_\mathrm{c} \in \child(Z_\mathrm{p}) \\ Z_\mathrm{c} = X \times I_\mathrm{c}}}
    \sum_{a_\mathrm{c} = 0}^{m_t}
    q^{(t)}_{a_\mathrm{c},a_\mathrm{p}}(I_\mathrm{c}, I_\mathrm{p})\, \mu_{a_\mathrm{c}, \vkappa}(Z_\mathrm{c}),
\end{equation}
with the coefficients $q^{(t)}_{a_\mathrm{c},a_\mathrm{p}}(I_\mathrm{c}, I_\mathrm{p})
  := L_{I_\mathrm{p},a_\mathrm{p}}(\xi_{I_\mathrm{c},a_\mathrm{c}}^{(m_t)})$.

\noindent \emph{Space-time M2M}:
If the children of a box $Z_\mathrm{p} = X_\mathrm{p} \times I_\mathrm{p}$ with $X_\mathrm{p} = X_\mathrm{p,1} \times X_\mathrm{p,2} \times X_\mathrm{p,3}$ are refined in space and time, there holds
\begin{equation} \label{eq:st_m2m}
  \mu_{a_\mathrm{p}, \vnu}(Z_\mathrm{p})
  = \sum_{\substack{Z_\mathrm{c} \in \child(Z_\mathrm{p}) \\ Z_\mathrm{c} = X_\mathrm{c} \times I_\mathrm{c}}}
    \sum_{a_\mathrm{c} = 0}^{m_t} \sum_{\vkappa \leq \vnu}
    q^{(t)}_{a_\mathrm{c},a_\mathrm{p}}(I_\mathrm{c}, I_\mathrm{p}) \,
    q^{(\vx)}_{\vkappa,\vnu}(X_\mathrm{c}, X_\mathrm{p}) \,
    \mu_{a_\mathrm{c}, \vkappa}(Z_\mathrm{c}).
\end{equation}
Here, the inequality $\vkappa \leq \vnu$ is understood componentwise, $q^{(t)}_{a_\mathrm{c},a_\mathrm{p}}(I_\mathrm{c}, I_\mathrm{p})$ are the coefficients appearing in~\eqref{eq:temp_m2m}, and ${q^{(\vx)}_{\vkappa,\vnu}(X_\mathrm{c}, X_\mathrm{p}) := \prod_j q^{(\vx)}_{\kappa_j,\nu_j}(X_{\mathrm{c},j},X_{\mathrm{p},j})}$ with
\begin{equation*}
  q^{(\vx)}_{\kappa_j,\nu_j}(X_{\mathrm{c},j},X_{\mathrm{p},j})
    := \frac{\lambda_{\kappa_j}}{m_\vx+1} \sum_{n=0}^{m_\vx}
      T_{X_{\mathrm{p},j},\nu_j} (\xi_{X_{\mathrm{c},j},n}^{(m_\vx)})
      T_{X_{\mathrm{c},j},\kappa_j} (\xi_{X_{\mathrm{c},j},n}^{(m_\vx)}),
\end{equation*}
where $\lambda_0 = 1$ and $\lambda_k=2$ for all $k\geq 1$.

% \noindent \emph{Temporal L2L}:
% If the child $Z_\mathrm{c} = X \times I_\mathrm{c}$ of a box $Z_\mathrm{p}= X \times I_\mathrm{p}$ is only refined in time, we compute
% \begin{equation} \label{eq:temp_l2l}
%   \widetilde{\lambda}_{b_\mathrm{c}, \vnu}(Z_\mathrm{c})
%   = \sum_{b_\mathrm{p}=0}^{m_t}
%     q^{(t)}_{b_\mathrm{c},b_\mathrm{p}}(I_\mathrm{c},I_\mathrm{p}) \lambda_{b_\mathrm{p},\vnu}(Z_\mathrm{p}),
% \end{equation}
% with $q^{(t)}_{b_\mathrm{c},b_\mathrm{p}}(I_\mathrm{c},I_\mathrm{p})$ as in \eqref{eq:s2m}. The result $\widetilde{\vlambda}(Z_\mathrm{c})$ is added to $\vlambda(Z_\mathrm{c})$.

% \noindent \emph{Space-time L2L}:
% For a child $Z_\mathrm{c} = X_\mathrm{c} \times I_\mathrm{c}$ of a box $Z_\mathrm{p} = X_\mathrm{p} \times I_\mathrm{p}$ we set
% \begin{equation} \label{eq:st_l2l}
%   \widetilde{\lambda}_{b_\mathrm{c}, \vnu_\mathrm{c}}(Z_\mathrm{c})
%   = \sum_{b_\mathrm{p}=0}^{m_t} \sum_{\substack{\vnu_\mathrm{p} \geq \vnu_\mathrm{c} \\ |\vnu_\mathrm{p}|\leq m_x}}
%     q^{(t)}_{b_\mathrm{c},b_\mathrm{p}}(I_\mathrm{c},I_\mathrm{p})
%     q^{(\vx)}_{\vnu_\mathrm{c},\vnu_\mathrm{p}}(X_\mathrm{c}, X_\mathrm{p})
%     \lambda_{b_\mathrm{p},\vnu_\mathrm{p}}(Z_\mathrm{p}),
% \end{equation}
% with the same coefficients $q^{(\vx)}_{\vnu_\mathrm{c},\vnu_\mathrm{p}}(X_\mathrm{c}, X_\mathrm{p})$ as for \eqref{eq:st_m2m}. As before, $\widetilde{\vlambda}(Z_\mathrm{c})$ is added to $\vlambda(Z_\mathrm{c})$.

The FMM presented in Algorithm~\ref{alg:fmm} reduces the runtime complexity of the matrix-vector multiplication~$\mat{V}_h \vw$ from~$\mathcal{O}((E_t\, E_\vx)^2)$ to $\mathcal{O}(m_t^2 m_\vx^4 E_t E_\vx)$ if the temporal and spatial mesh sizes~$h_t$ and $h_x$ of $\Sigma_h$ satisfy~$h_t \sim h_x^2$, see e.g.,~\cite[Section~5.4]{Tausch2007}. In~\cite{MessnerSchanzTausch14} an additional nearfield compression is provided for meshes whose temporal mesh sizes are too large. Such a nearfield compression is not considered in this work.

\begin{algorithm}[H]
  \caption{Space-time FMM for the approximate evaluation of $\vf = \mat{V}_h \vw$} \label{alg:fmm}
  \begin{algorithmic}[1]
  \State Choose the parameters $n_{\mathrm{max}}$, $c_\mathrm{st}$, $n_\mathrm{tr}$ and the expansion degrees $m_t$ and $m_\vx$.
  \State Construct the box cluster tree $\mathcal{T}_\Sigma$ and determine the sets $\mathcal{N}(Z)$ and $\mathcal{I}(Z)$ for all $Z \in \mathcal{T}_\Sigma$ according to \eqref{eq:def_st_nearfield} and \eqref{eq:def_st_interaction_list}.
  \State Initialize $\vf = \vzero$.
  \State \Comment{Forward transformation}
  \For{all leaves $Z \in \mathcal{L}_\Sigma$}
    \State S2M: Compute $\vmu(Z)$ by \eqref{eq:s2m}.
  \EndFor

  \For{all levels $\ell = p(\mathcal{T}_\Sigma) - 1$, \ldots, $2$}
    \For{all non-leaf boxes $Z_\mathrm{p} \in \mathcal{T}_\Sigma^{(\ell)}$}
      \If{children of $Z_\mathrm{p}$ are refined only in time}
        \State Temporal M2M: Compute $\vmu(Z_\mathrm{p})$ by \eqref{eq:temp_m2m}.
      \Else
        \State Space-time M2M: Compute $\vmu(Z_\mathrm{p})$ by \eqref{eq:st_m2m}.
      \EndIf
    \EndFor
  \EndFor

  \State \Comment{Multiplication phase}
  \For{all boxes $Z_{\mathrm{tar}} \in \mathcal{T}_\Sigma$}
    \State Initialize $\vlambda(Z_{\mathrm{tar}})=\vzero$.
    \For{all boxes $Z_{\mathrm{src}} \in \mathcal{I}(Z_\mathrm{tar})$}
      \State M2L: Update $\vlambda(Z_{\mathrm{tar}})$ by adding the result from \eqref{eq:m2l}.
    \EndFor
  \EndFor

  \State \Comment{Backward transformation}
  \For{all levels $\ell = 3$, \ldots, $p(\mathcal{T}_\Sigma)$}
    \For{all boxes $Z_\mathrm{c} \in \mathcal{T}_\Sigma^{(\ell)}$}
      \If{$Z_\mathrm{c}$ results from its parent $Z_\mathrm{p}$ by a purely temporal refinement}
        \State Temporal L2L: Update $\vlambda(Z_\mathrm{c})$ using $\vlambda(Z_\mathrm{p})$.
      \Else
        \State Space-time L2L: Update $\vlambda(Z_\mathrm{c})$ using $\vlambda(Z_\mathrm{p})$.
      \EndIf
    \EndFor
  \EndFor

  \For{all leaves $Z \in \mathcal{L}_\Sigma$}
    \For{all $k_t$ and $k_\vx$ such that $\sigma_{k_t,k_\vx} \in Z$}
      \State L2T: Update $f_{k_t,k_\vx}$ by adding the result from \eqref{eq:l2t}.
    \EndFor
  \EndFor

  \State \Comment{Nearfield evaluation}
  \For{all leaves $Z_{\mathrm{tar}} \in \mathcal{L}_\Sigma$}
    \For{all $Z_{\mathrm{src}}$ in the nearfield $\mathcal{N}(Z_{\mathrm{tar}})$}
      \State Update $f|_{\hat{Z}_\mathrm{tar}} \mathrel{+}=
        \mat{V}_h|_{\hat{Z}_{\mathrm{tar}} \times \hat{Z}_{\mathrm{src}}} \vw|_{\hat{Z}_{\mathrm{src}}}$.
    \EndFor
  \EndFor
  \end{algorithmic}
  \end{algorithm}

\section{Parallelization of the space-time FMM in shared and distributed memory}
\label{sec:04_parallelization}
We present a novel parallel implementation of the space-time FMM presented in Section~\ref{sec:03_stFMM}.
Our distributed parallelization strategy relies on a decomposition of the one-dimensional temporal tree~$\mathcal{T}_I$, see Section~\ref{subsec:nfia}, into locally essential trees~(LET) \cite{WarrenSalmon1993}, which are distributed among the available MPI processes. A LET is a local part of a tree extended by fractions of the global tree which are required in the collaboration of the processes. The FMM operations are assigned to the processes according to these temporal subtrees. The necessary inter-process communication is also handled clusterwise in the temporal tree in a one-directional way mainly between successive processes due to the causality of the boundary integral operators. This clusterwise communication leads to a small number of communication events of reasonable size.

A distributed space-time cluster tree is assembled collaboratively based on the temporal tree. Certain parts of the FMM can be computed locally on the space-time subtrees, while other parts of the computations depend on the results from other subtrees, i.e.,~processes. We employ a small scheduling system to keep track of the dependencies.
Relations between the temporal clusters within the global distributed tree define dependencies which are used to decompose the computation in Algorithm~\ref{alg:fmm} into tasks which can be executed asynchronously. During the matrix vector multiplication a scheduler goes through a list of available tasks and executes those with fulfilled dependencies. This enables an asynchronous parallelization instead of bulk-synchronous approach which is often used to parallelize the FMM. Let us describe the individual phases of the computation in more details.

\subsection{Preprocessing}

\begin{figure}[t]
  \begin{center}
      \includegraphics[width=1\textwidth]{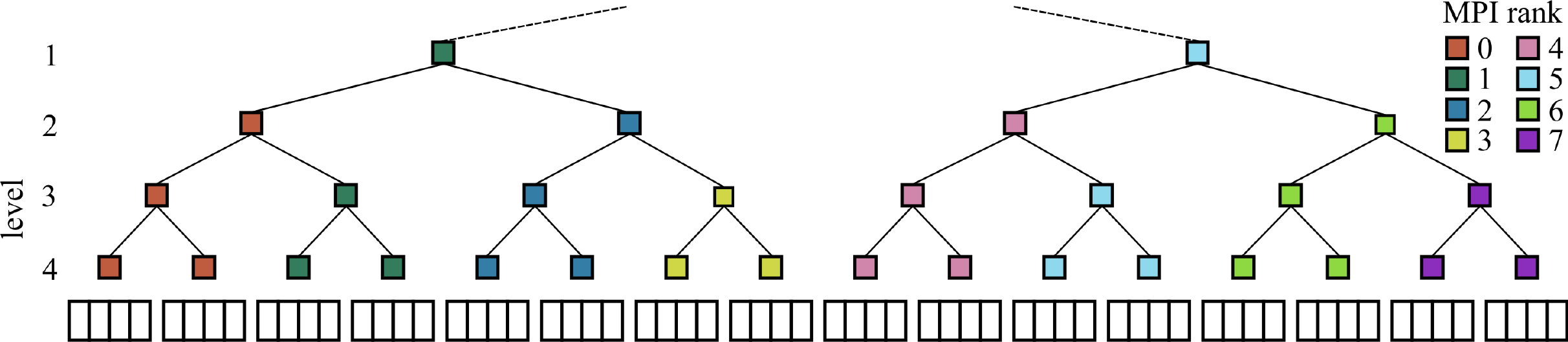}
  \end{center}
  \caption{16 initial time-slices each containing 4 time-steps. The distribution tree is built above these time-slices and split among available MPI processes.}
  \label{fig:slices_tree}
\end{figure}

The preprocessing phase consists of the assembly of the temporal cluster tree and its decomposition, the assembly of the distributed space-time mesh, and the creation of the distributed space-time cluster tree. We start by splitting the global time interval $(0,T)$ and the related time-steps into time-slices. Each of the slices contains multiple time-steps (see Figure~\ref{fig:slices_tree}). Next, a binary tree is built by a recursive bisection of the global interval and the time-slices are assigned to the nodes of the tree accordingly. We distribute the individual nodes of the temporal cluster tree among the available $N_{\mathrm{proc}}$ MPI processes. This temporal decomposition will drive the distribution of the space-time cluster tree later on and will determine the required communication between processes. In particular, all processes will execute the FMM operations related to space-time clusters that are associated with temporal clusters, for which they are responsible. Multiple distribution approaches can be employed; in our experiments we try to reduce inter-process communication and to obtain a reasonable load balance by determining the distribution starting from the finest tree level according to the following strategy:
\begin{itemize}
    \item On level $\lceil\log_2{N_{\mathrm{proc}}}\rceil$ and finer levels of the temporal cluster tree there are more clusters than processes. We distribute the clusters and the related time-slices among all $N_{\mathrm{proc}}$ processes in ascending order and as uniformly as possible.
    \item Level $\lceil\log_2{N_{\mathrm{proc}}}\rceil - 1$ is the first level where the number of clusters is less than the number of available processes. Here we assign a cluster to the process that handles its left child. In this way we improve the load balancing because the temporal interaction list of the left child is smaller than the one of the right child, see~\eqref{eq:03:ILtemp}.
    \item On each level $\ell$ with $\ell < \lceil\log_2{N_{\mathrm{proc}}\rceil} - 1$ there are again less clusters than processes, so we split the processes into $2^\ell$ groups of ascending order to assign the $2^\ell$ cluster of the level. From each group we pick a process, which has been responsible for the smallest number of clusters so far, and assign it to the related cluster. Again we aim at improving the load balancing.
\end{itemize}

Figure~\ref{fig:slices_tree} gives an example of a distribution based on this strategy. Please check the coloring of the nodes for the assignment of processes. On levels 3 and 4 the process assignment is driven by the first rule and clusters on these levels are equally distributed among the available MPI processes. On level 2 we employ the second rule and assign each cluster to the MPI process handling its left child to improve the load balancing. Finally, clusters on level 1 are distributed according to the third rule. Processes are split into two groups (0 -- 3 and 4 -- 7) and from each group we select a process with the smallest number of clusters. In this example the selected processes are 1 and 5.

\begin{figure}[b]
  \begin{center}
      \includegraphics[width=0.6\textwidth]{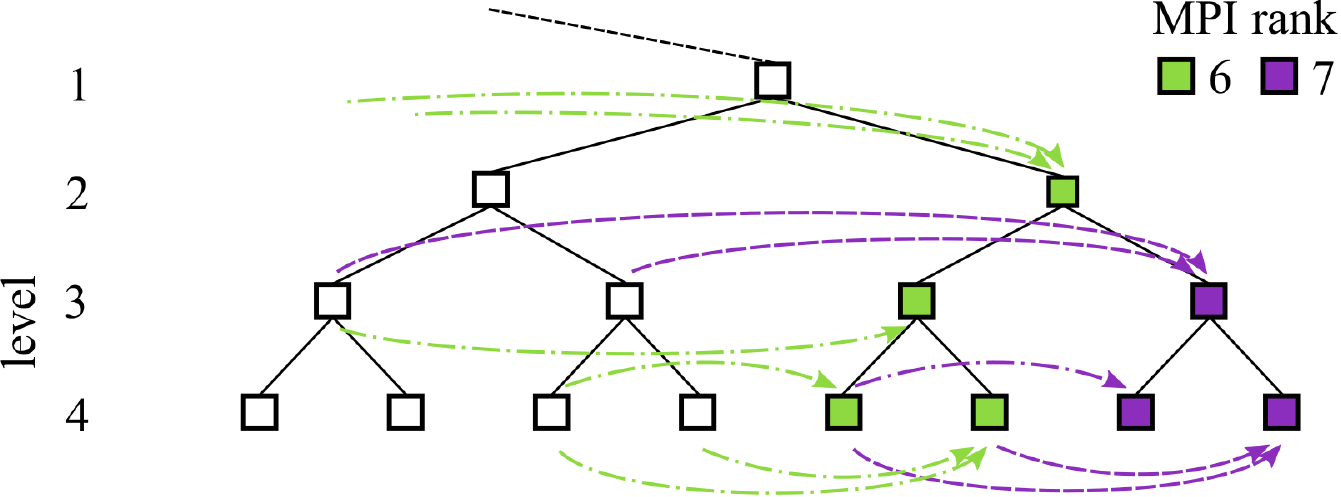}
  \end{center}
  \caption{M2L operations on MPI processes with ranks 6 (green dash-dotted line) and~7 (purple dashed line) in the tree from Figure~\ref{fig:slices_tree}. For the sake of clarity, only the relevant part of that tree is depicted.}
  \label{fig:M2L_comm}
\end{figure}

Let us discuss the load balance of the example from Figure~\ref{fig:slices_tree} to motivate the suggested distribution strategy.
We can restrict the discussion to the M2L operations, since nearfield and M2L operations generate most computational effort of the FMM, and the nearfield operations are evenly distributed among the processes by the original distribution on the finest level.
We compare the efforts of the processes with ranks 6 and 7. According to the arrows in Figure~\ref{fig:M2L_comm} indicating the M2L operations, the related total efforts $W_6$ and $W_7$ are the sums of the individual efforts $M2L(\ell)$ on the different temporal levels $\ell$
\begin{align*}
  W_6 = 3 \cdot M2L(4) + 1 \cdot M2L(3) + 2 \cdot M2L(2),\quad
  W_7 = 3 \cdot M2L(4) + 2 \cdot M2L(3).
\end{align*}
Note the imbalance of one and two M2L operations on level 3 due to the differently sized interaction lists of the clusters. We normalize the effort on level 2 by setting $M2L(2)=1$.
The effort of M2L operations on levels 3 and 4 depends on the underlying space-time cluster tree. Exemplarily, we discuss the two most relevant scenarios.
If the cluster refinement from level 2 to~3 is purely temporal and there is a space-time refinement from level 3 to 4, we have the efforts $M2L(3)=1$ and $M2L(4)=4$, where we assume that each cluster has $4$ spatial children, because the spatial parts of the clusters resolve the spatial surface, not the volume. Then the total efforts are comparable as
\begin{align*}
  W_6 = 3 \cdot 4 + 1 \cdot 1 + 2 \cdot 1 = 15, \quad
  W_7 = 3 \cdot 4 + 2 \cdot 1 = 14.
\end{align*}
If there are a space-time cluster refinement from level 2 to~3 and a pure temporal refinement from level 3 to 4, we have the efforts $M2L(3)=4$ and $M2L(4)=4$ and again the total efforts are comparable as
\begin{align*}
  W_6 = 3 \cdot 4 + 1 \cdot 4 + 2 \cdot 1 = 18, \quad
  W_7 = 3 \cdot 4 + 2 \cdot 4 = 20.
\end{align*}
Let us mention that the process with rank 0 does not have to execute any M2L operation in this simple example. This imbalance does not really matter for a larger number of processes. Finally, it is important to note that the cost of M2L operations increases exponentially with the number of spatial refinements, i.e.~with increasing level. In particular, the total effort of M2L operations is dominated by the effort on the fine levels. Hence, small differences in the upper part of the trees as discussed above are insignificant, if there are additional fine levels in the cluster tree which are evenly distributed.

\begin{figure}
  \begin{center}
      \includegraphics[width=1\textwidth]{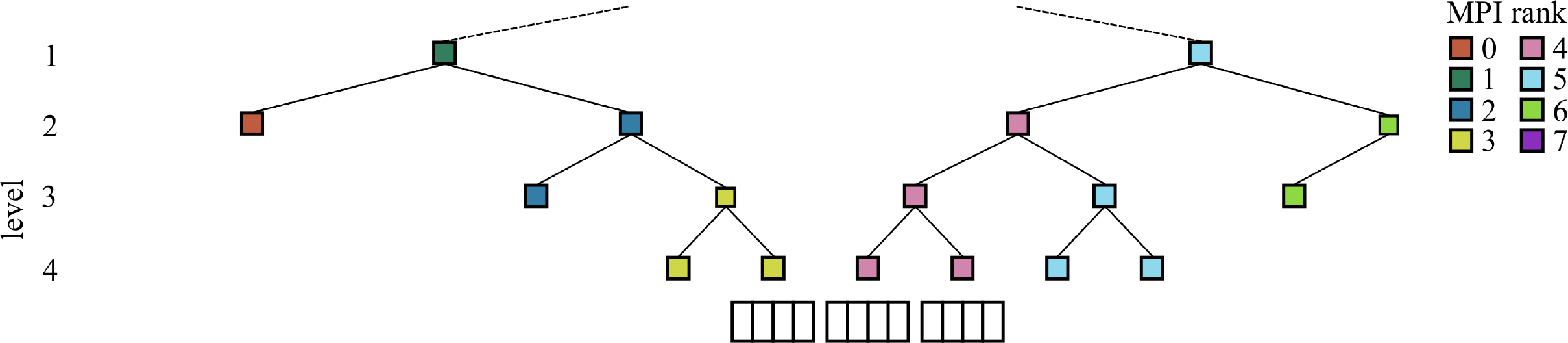}
  \end{center}
  \caption{Locally essential tree of the MPI process with rank 4 from the example in Figure~\ref{fig:slices_tree}. The process holds all the data associated with its nodes of the tree (denoted by the pink color). The remaining nodes of its LET contain only necessary information, e.g.~the ranks of the remote processes that are responsible for the data. The drawn time-steps are those for which process 4 stores the related part of the space-time mesh. Note that we include some additional time-steps which are contained in a cluster assigned to the process with MPI rank 3 for a more efficient nearfield computation.}
  \label{fig:LET}
\end{figure}

The assignment of temporal clusters to MPI processes determines the part of the temporal tree which is relevant for a process~$p$. Besides the assigned clusters, process~$p$ also needs to know about clusters that are needed for nearfield computations and clusters which are related to FMM operations requiring communication from or to process~$p$. All these clusters form the LET of a process~$p$ (see Figure~\ref{fig:LET} for the LET of the MPI process with rank 4 from the example in Figure~\ref{fig:slices_tree}). The clusters in the LET for which the process is not directly responsible contain the ranks of the responsible remote processes. Since the temporal tree defines the distribution of work among processes we will call it the \emph{scheduling tree}.
Next, we combine the time-steps in the individual time-slices with the spatial mesh to create a distributed space-time tensor product mesh. Each process gets the local part of the full space-time mesh related to time-steps contained in leaf clusters of the scheduling tree for which it is directly responsible as well as time-steps in the temporal nearfield of its clusters. The latter are included to reduce communication during the assembly of the nearfield blocks of the matrix.

A distributed 4D space-time box cluster tree is created level-wise top down as described in Section~\ref{sec:4d:tree} and in Algorithm~\ref{alg:box_cluster_tree}. The nodes of the tree are assigned to the processes based on the assignment of their temporal components defined in the scheduling tree. Each node of the scheduling tree stores the information about the corresponding space-time clusters. The depth of the space-time tree may be larger than the depth of the temporal tree. In such a case the temporal tree is locally extended accordingly. When building the upper part of the space-time cluster tree, where clusters contain elements from local meshes assigned to multiple processes, communication is required but can be limited to a reasonable amount. However, the construction of the lower parts of the tree is done independently on each process and just some synchronization with direct neighbors is carried out to set up the communication.

\subsection{Matrix-vector multiplication}

Our parallelization approach in shared and distributed memory is based on a data driven model, instead of a bulk-syn\-chro\-nous parallelization often used in scientific codes. In \cite{abduljabbar2014asynchronous} such an approach is used for the parallelization of an FMM for particle simulations using \texttt{Charm++}. We use a different strategy for work distribution and execution and a custom scheduler in combination with OpenMP \texttt{tasks} in order to avoid dependencies on external software. Our approach has similarities to \cite{Aguetal2014} which we will discuss at the end of the current section. To distinguish our own tasks from  OpenMP \texttt{tasks}, the latter ones will be denoted in monospaced font from now on.

Our top-level scheduler is based on the FMM operations and dependencies with respect to the temporal scheduling tree. At this stage just consider Algorithm~\ref{alg:fmm} reduced to the scheduling tree~$\mathcal{T}_I$ and the related temporal operations. If such a temporal FMM operation is called, OpenMP \texttt{tasks} of all related space-time operations of the attached space-time clusters are generated to implement shared memory parallelization and to control the granularity of the \texttt{tasks}.

The strategy of the temporal scheduler is as follows:
Most tasks are related to pairs of clusters.
Certain temporal tasks do not have dependencies and can be executed at any time, while some temporal tasks depend on the results of other temporal tasks. We distinguish local and remote dependencies. If a computation is done on the same process, the dependency is resolved locally.
If a tasks requires data from a computation of another process, the calculation can start at earliest after the communication has taken place.
The temporal scheduler prioritizes those tasks which other tasks depend on. In case none of these tasks can be executed at some time, we avoid idle times by scheduling independent tasks.

In more details, we start by decomposing the temporal variant of Algorithm~\ref{alg:fmm} into tasks and defining their mutual dependencies with respect to the temporal scheduling tree. We distinguish the following lists of tasks:
\begin{itemize}
    \item \emph{M-list} -- S2M and M2M operations including send operations of the computed moments to the parent and the clusters in the interaction list,
    \item \emph{M2L-list} -- transformations of moments into local contributions (M2L) including possible downward send operations or evaluations of the local contributions in case of a leaf (L2T),
    \item \emph{L-list} -- translations of local contributions from parent (L2L) including possible downward send operations or evaluations of the local contributions in case of a leaf (L2T),
    \item \emph{N-list} -- execution of the nearfield operations.
\end{itemize}
The tasks in a list are ordered in a way which is advantageous for our parallel execution.
Within each list a task can be identified from the related cluster. Thus we just store temporal clusters in these lists.
Each cluster is assigned to one or more of the lists depending on its position in the tree. This enables us to specify dependencies for individual clusters (tasks):
\begin{itemize}
    \item Tasks of non-leaf clusters in the M-list depend on the completion of the M-list operations of their children.
    \item Tasks of clusters in the M2L-list depend on the completion of M-list operations of clusters in their interaction lists (see~\eqref{eq:03:ILtemp}).
    \item Tasks of clusters in the L-list depend on the completion of their parents' M2L- and L-list operations.
    \item Nearfield tasks of clusters in the N-list have no dependencies.
\end{itemize}

\begin{algorithm}[t]
  \caption{Parallel space-time FMM for the approximate evaluation of $\vf = \mat{V}_h \vw$} \label{alg:par_fmm}
  \begin{algorithmic}[1]
      \State Fill the \texttt{M\_list}, \texttt{M2L\_list}, \texttt{L\_list}, and \texttt{N\_list}; Initialize $\vf = \vzero$.

      \Pragma[OpenMP \texttt{parallel} region]
      \Pragma[OpenMP \texttt{single} section]
      \State \Call{StartMPIReceiveOperations}{ }
      \While{the lists are not empty}
        \State \Call{CheckMPIForReceivedData}{ }
        \State \texttt{[cluster, list]}
        \makeatletter
        \Statex \hskip\ALG@thistlm $\qquad$= \Call{FindNextCluster}{\texttt{M\_list}, \texttt{L\_list}, \texttt{M2L\_list}, \texttt{N\_list}}
        \makeatother

        \If{\texttt{n\_generated\_tasks > threshold}}
          \State Suspend the execution of scheduling task using the \texttt{taskyield} construct \label{alg:line:taskyield}
        \EndIf

        \If{\texttt{list} == 0}
          \State \Call{CreateOpenMPMListTask}{\texttt{cluster}}
          \State \Call{RemoveClusterFromList}{\texttt{cluster}, \texttt{M\_list}}
        \ElsIf{\texttt{list} == 1}
          \State \Call{CreateOpenMPLListTask}{\texttt{cluster}}
          \State \Call{RemoveClusterFromList}{\texttt{cluster}, \texttt{L\_list}}
        \ElsIf{\texttt{list} == 2}
          \State \Call{CreateOpenMPM2LListTask}{\texttt{cluster}}
          \State \Call{RemoveClusterFromList}{\texttt{cluster}, \texttt{M2L\_list}}
        \ElsIf{\texttt{list} == 3}
          \State \Call{CreateOpenMPNListTask}{\texttt{cluster}}
          \State \Call{RemoveClusterFromList}{\texttt{cluster}, \texttt{N\_list}}
        \EndIf

      \EndWhile
      \EndPragma
      \EndPragma

  \end{algorithmic}
\end{algorithm}

Note that each process creates only the parts of these lists which are relevant to its local part of the scheduling tree and to its LET, respectively. We do not describe this restriction explicitly but it results from the local scheduling tree naturally.
A simplified parallel matrix-vector multiplication algorithm is described in Algorithm~\ref{alg:par_fmm}. We start by filling the above mentioned lists by clusters of the scheduling tree according to the presented rules. In addition to the distributed parallelization with MPI, we make use of OpenMP thread parallelization within each process. First, we enter the OpenMP \texttt{parallel} region and create a \texttt{single} section to ensure that only one thread will execute the main scheduling loop. The routine \textsc{StartMPIReceiveOperations()} creates a non-blocking receive operation using the \texttt{MPI\_Irecv()} function for every temporal cluster in the lists requiring data from remote processes. Individual clusters and operations are distinguished using the \texttt{tag} argument of the MPI function.

At the beginning of each iteration in the while loop of Algorithm~\ref{alg:par_fmm} the scheduling thread calls \textsc{CheckMPIForReceivedData()} which uses \texttt{MPI\_Testsome()} to check for new data received from remote processes and updates the dependencies of the respective tasks and temporal clusters. Then the scheduling thread iterates through the lists calling \textsc{FindNextCluster} to find the next cluster or rather task ready to be executed, i.e.,~one with all dependencies fulfilled. Due to the succession of dependencies, first the M-list tasks are checked, then the L-list tasks, followed by the M2L-list tasks. Since the nearfield tasks are independent of all other tasks, the N-list is traversed last. If the scheduling thread finds a task, it creates an OpenMP \texttt{task} for executing the corresponding operations on the cluster and removes the cluster from the list using \textsc{RemoveClusterFromList}.
Here the OpenMP task scheduler serves as a buffer for the initialized \texttt{tasks} and allows for a finer granularity of the tasks and a more efficient parallelization.

In order to avoid collisions of the OpenMP \texttt{tasks} during the memory access, the generated OpenMP \texttt{task} is subject to additional dependencies specified using the OpenMP \texttt{depend} clause. E.g., a task generated from the M-list depends on all previously generated \texttt{tasks} where the cluster has the same parent as the current cluster (to avoid collisions during the M2M operations) and the \texttt{tasks} where the cluster is in the interaction lists of the same cluster as the current cluster (to avoid collisions during the M2L operations).

When executing a task from one of the four lists, the related operations in the space-time cluster tree are processed. To improve the granularity of the shared memory parallelization, additional finer level OpenMP \texttt{tasks} are created by iterating through the associated space-time clusters using the \texttt{taskloop} construct. After these associated finer level OpenMP \texttt{tasks} are completed, the dependencies in the temporal scheduler are updated. If the depending cluster is owned by a remote process, the required data are sent using the non-blocking \texttt{MPI\_Isend()} operation.

After all lists are empty and all OpenMP \texttt{tasks} are completed, every process has computed its local part of the matrix-vector product $\vf = \mat{V}_h \vw$.

Notice that the scheduling thread may suspend the execution of the \texttt{while} loop if the number of already generated tasks is greater than a certain threshold and join other threads in executing the generated \texttt{tasks} (see line \ref{alg:line:taskyield} in Algorithm~\ref{alg:par_fmm}). However, the \texttt{taskyield} construct is a non-binding request and may result in a \texttt{no} \texttt{operation}. In this context, a suspension of the scheduling \text{task} is not observed, e.g., in the GCC compiler (v9.3). Therefore, we mainly focus on the Intel compiler in Section~\ref{sec:05_experiments} on our numerical experiments.

\paragraph{Note} One could choose another task-based approach to hybrid OpenMP-MPI parallelization. Instead of providing a custom scheduler run by a dedicated thread, all tasks could be created at once including tasks responsible for data receiving and sending. The correct execution order would be ensured using the OpenMP \texttt{depend} clause. Unfortunately, due to the above-described characteristics of the \texttt{taskyield} construct and a limited number of OpenMP tasks we cannot prevent deadlocks when only receive or send operations are posted on the processes. Alternatively, one could rely on special compilers supporting task suspending such as OmpSS \cite{SchuchartEtAl}. However, we decided to implement our own simple scheduler to reduce the number of external dependencies. This also enables us to better control the granularity of computation.

\paragraph{Note}
Similar to \cite{Aguetal2014} our parallelization relies on a task based reformulation of the FMM.
While~\cite{Aguetal2014} describes a shared memory parallelization and only mentions a possible extension to MPI, our approach is designed for a native distributed memory parallelization. For that purpose we use two levels of tasks.
Our top level tasks realize the distributed parallelization and related one-directional communication based on the temporal cluster tree. They additionally cover all dependencies of FMM operations.
The top level tasks group the space-time FMM operations according to the underlying temporal tree which allows to exploit the temporal structure in the communication and the FMM operations.
The grouping induces a coarse granularity which we refine by our second level OpenMP tasks.
Instead, in \cite{Aguetal2014} a single task scheduler (StarPU~\cite{StarPU}) is used, which is more advanced than OpenMP and handles all dependencies, and the granularity of tasks is adjusted by collecting clusters into larger blocks. The approach from \cite{Aguetal2014} was extended in \cite{agullo-report} where the authors compare several approaches for task-based distributed memory parallelization of FMM using StarPU as a scheduler. In most of the presented approaches, the inter-node communication is fully delegated to the StarPU runtime system. While their runtime system automatically derives the necessary communication from the spatial decomposition and dependencies between individual operations, we tailor our communication pattern to the temporal component and benefit from simpler structures and the one-sided communication. This way we can use a relatively simple and lightweight scheduler without dependencies on external software.

\paragraph{Note} Since we based our MPI parallelization on a decomposition of the temporal scheduling tree, the distributed memory parallelization is limited by the total number of time-steps. Extending the distributed parallelization into a spatial dimension would thus improve the scalability on large machines. The simplest approach that does not significantly disrupt the existing implementation would be to start from an existing decomposition and to replace each MPI process by a group of processes. Within this group, the processes would collaborate on the local space-time operations and a certain process from the group would be responsible for the communication with other groups using the same communication structure as the current algorithm. This hierarchical nature of the communication could exploit the topology of the computer cluster and thus further reduce the communication time.

\section{Numerical experiments}
\label{sec:05_experiments}
The parallel space-time FMM algorithm from Section~\ref{sec:04_parallelization} has been implemented in the publicly available \cpp{} library besthea~\cite{besthea}.
As we aim for space-time adaptive methods, we do not exploit the block-Toeplitz structure of the global BEM matrices which is only present for uniform meshes. In the current state, we do not apply a compression of the temporal nearfield blocks as developed  in \cite{MessnerSchanzTausch14}, which may improve the performance.
To evaluate the efficiency of the presented space-time FMM we carried out numerical experiments using the Salomon and Barbora clusters at IT4Innovations National Supercomputing Center in Ostrava, Czech Republic. The Salomon cluster consists of 1009 compute nodes equipped with two 12-core Intel Xeon E5-2680v3 processors and 128 GB of RAM. The theoretical peak performance of the cluster is 2 PFLOP/s. The Barbora cluster consists of 201 computational nodes equipped with two 18-core Intel Cascade Lake 6240 CPUs and 192 GB of~RAM. Nodes within both clusters are interconnected using the InfiniBand network. On Salomon we used the Intel Compiler v19.1.1 and the Intel Math Kernel Library (MKL) v2020.1 unless stated otherwise, while on Barbora we employed the Intel Compiler v19.1.3 and MKL v2020.4. The affinity of the threads to cores was set using the variable \texttt{KMP\_AFFINITY=granularity=core,compact}, which guarantees that the threads will stay within a single socket of the two-socket system when possible. The details of how to reproduce the numerical results are provided in the text below and in the software repository \cite{besthea}.

\subsection{Shared memory performance}

The performance of our OpenMP parallelization of the space-time FMM was tested on up to 36 cores on a node of the Barbora cluster. The scalability was tested for the Dirichlet problem and the boundary element method described in Section~\ref{sec:02_bie_bem}. As the Dirichlet datum we use
\begin{equation}
    u(\vx, t) = G_{\alpha}(\vx-\vy^*,t) \quad \text{for } (\vx,t)\in \Sigma,
    \label{eq:dirdata}
\end{equation}
where $\vy* = (1.5, 1.5, 1.5)$ and $\alpha=1$.
The lateral boundary of the space-time domain $(-0.5, 0.5)^3\times (0, 0.25)$ was discretized into 1\,536 spatial and 64 temporal elements, resulting in 98\,304 space-time boundary elements in total. The temporal elements were distributed among 16 time-slices. The space-time leaf cluster size was limited using the value $n_{\mathrm{max}}=80$ (see Algorithm~\ref{alg:box_cluster_tree}) and the orders of the Lagrange and Chebyshev polynomials were both set to $m_t=m_\vx = 6$. Finally, we used $c_{\mathrm{st}}=0.9$ in the relation~\eqref{eq:rel_h_t_h_x} controlling the spatial and temporal box sizes and the cutting parameter $n_{\mathrm{tr}}=5$ in~\eqref{eq:03:Iarea}. Note that here and in all other experiments we have chosen the parameters for the FMM such that they do not affect the approximation quality of the BEM (see, e.g.,~\cite{costabel:1990,DohNiiSte2019} for results regarding the approximation quality) and used GMRES without a preconditioner to solve the considered linear systems of equations.

The assembly times of the single and double layer matrices, respectively, and the time per GMRES iteration are given in Table~\ref{tab:threading} with respect to the number of threads. Up to 18 threads on a single socket the efficiency of the matrix assembly is above 90\%, while for 36 threads spanning over two sockets the efficiency drops mainly for the assembly of the double layer matrix due to the less efficient memory access. The iterative solution, which is mostly composed of the matrix-vector multiplication, scales almost optimally even for both sockets fully occupied. The GMRES solver required 32 iterations to attain the relative precision~$10^{-8}$. The total computational time was reduced from 8282 s on one thread to 245 s on 36 threads.

% \begin{table}[t]
%     \begin{center}
%         \begin{tabular}{rrrrrrrr}
%             \toprule
%             & \# threads & 1 & 2 & 4 & 8 & 18 & 36 \\
%             \midrule
%             \multirow{2}{*}{$\matV_{h}$} & time $[s]$ & 509.9 & 263.6 & 138.9 & 69.9 & 32.9 & 16.6  \\
%            & efficiency [\%]& 100.0 & 96.7 & 91.8 & 91.1 & 86.1 & 85.5 \\
%            \midrule
%           \multirow{2}{*}{$\matK_{h}$} & time $[s]$ & 501.3 & 266.7 & 140.1 & 70.5 & 33.9 & 20.3\\
%            & efficiency [\%] & 100 & 94.0 & 89.5 & 88.9 & 82.1 & 68.6 \\
%            \midrule
%           \multirow{2}{*}{iteration} & time $[s]$& 225.8 & 113.2 & 56.7 & 28.4 & 12.8 & 6.6 \\
%           & efficiency [\%] & 100 & 99.8 & 99.6 & 99.5 & 97.9 & 95.4  \\
%           \bottomrule
%         \end{tabular}
%     \end{center}
%     \caption{Assembly times and time per GMRES iteration for different numbers of OpenMP threads and a problem with 98\,304 space-time surface elements (64 time-steps, 1\,536 spatial elements) using a single node of the Barbora cluster.}
%     \label{tab:threading}
% \end{table}

\begin{figure}[ht]
    \begin{center}
        \includegraphics[width=0.32\textwidth]{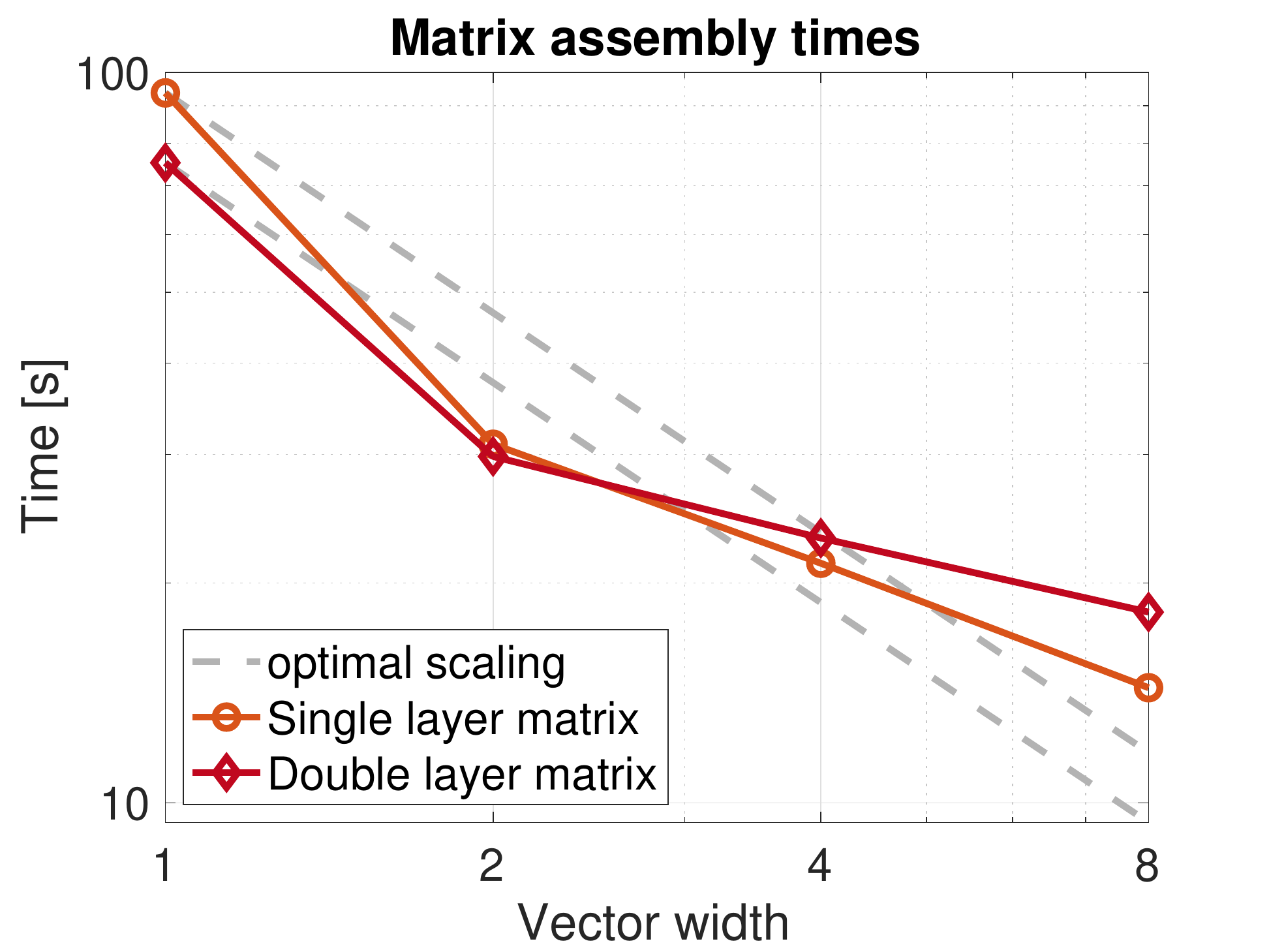}
        \includegraphics[width=0.32\textwidth]{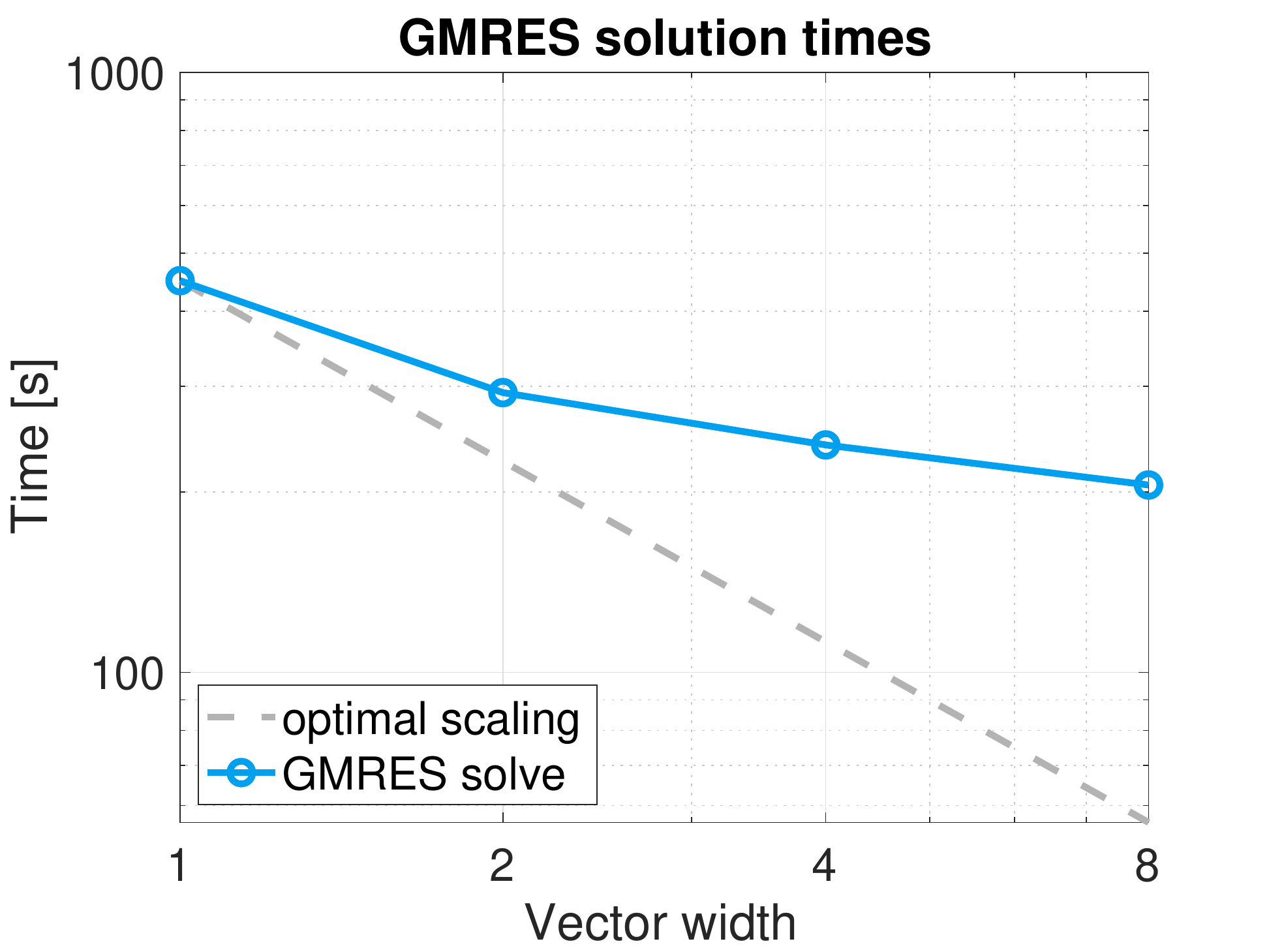}
        \includegraphics[width=0.32\textwidth]{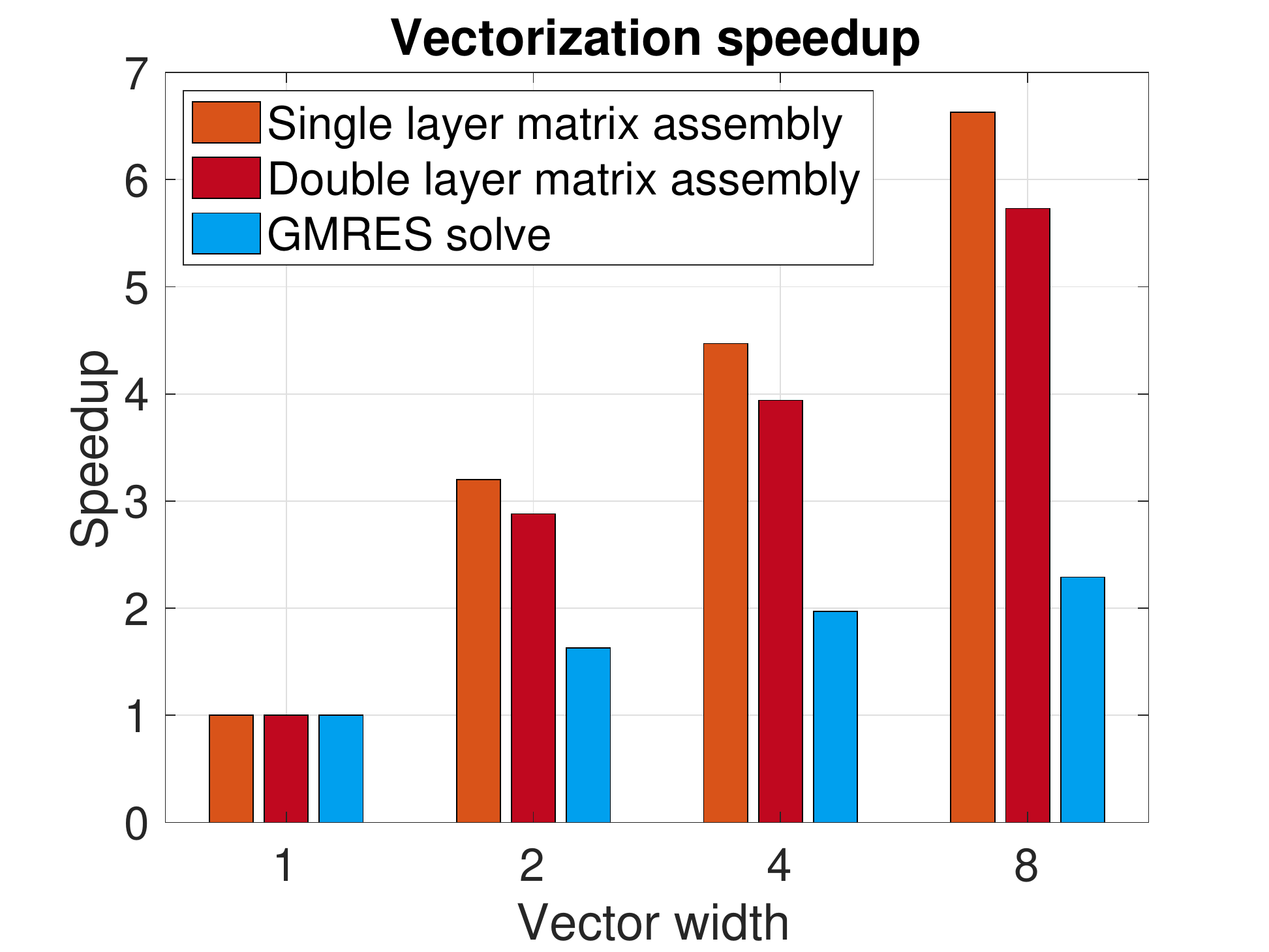}
    \end{center}
    \caption{Efficiency of the vectorization with respect to the width of the vector registers. The non-vectorized version (vector width equals one) was compiled using the compiler flags \texttt{-no-vec -no-simd -qno-openmp-simd}. The remaining versions were compiled using the flags \texttt{-xcore-avx512 -qopt-zmm-usage=high} with the vector width set by the \texttt{simdlen} OpenMP clause.}
    \label{fig:simd}
\end{figure}

On modern CPUs, the performance of a code highly depends on its ability to exploit SIMD vectorization. In our code we employ vectorization using OpenMP \texttt{simd} pragmas for the assembly of the nearfield part of the matrices, see \cite{zapletal2021semianalytic} for the details, and the M2L operations in \eqref{eq:m2l}. For the latter we use the efficient realization discussed in \cite[Section 4.3]{TauWec2009} and apply SIMD vectorization when computing the coefficients \eqref{eq:exp_coeffs_3d} or rather \eqref{eq:exp_coeffs_1d} and the occurring one-dimensional transforms. In addition, we use Intel MKL/BLAS to compute the matrix-vector products in the nearfield operations in \eqref{eq:eval_sl_nf_block}. The effect of the vectorization is demonstrated in Figure~\ref{fig:simd} where we present the scalability of the assembly of the system matrices (mainly involving the assembly of nearfield parts) and the iterative solution using GMRES with respect to the SIMD registers vector width. We used the same settings as in the previous experiment. Compared to the matrix assembly, the GMRES vector scalability is limited due to the complexity of the space-time FMM code for the matrix vector multiplication. While further optimization would probably be possible, the current code is a trade-off between the performance and readability. Nevertheless, the achieved speedup still reduces the solution time significantly.

\begin{table}[t]
    \begin{center}
        \begin{tabular}{rrrrrrrr}
            \toprule
            & \# threads & 1 & 2 & 4 & 8 & 18 & 36 \\
            \midrule
            \multirow{2}{*}{$\matV_{h}$} & time $[s]$ & 505.7 & 254.6 & 127.4 & 63.8 & 28.4 & 14.4  \\
           & efficiency [\%]& 100.0 & 99.3 & 99.3 & 99.2 & 98.8 & 97.6 \\
           \midrule
          \multirow{2}{*}{$\matK_{h}$} & time $[s]$ & 488.2 & 258.9 & 129.9 & 64.9 & 29.4 & 18.3\\
           & efficiency [\%] & 100.0 & 94.3 & 94.0 & 94.0 & 92.3 & 74.1 \\
           \midrule
          \multirow{2}{*}{iteration} & time $[s]$& 220.9 & 111.5 & 55.5 & 27.9 & 12.5 & 6.4 \\
          & efficiency [\%] & 100.0 & 99.1 & 99.5 & 98.8 & 98.3 & 95.8  \\
          \bottomrule
        \end{tabular}
    \end{center}
    \caption{Assembly times and time per GMRES iteration for different numbers of OpenMP threads and a problem with 98\,304 space-time surface elements (64 time-steps, 1\,536 spatial elements) using a single node of the Barbora cluster.}
    \label{tab:threading}
\end{table}

\begin{figure}[t]
    \begin{center}
        \includegraphics[width=1\textwidth]{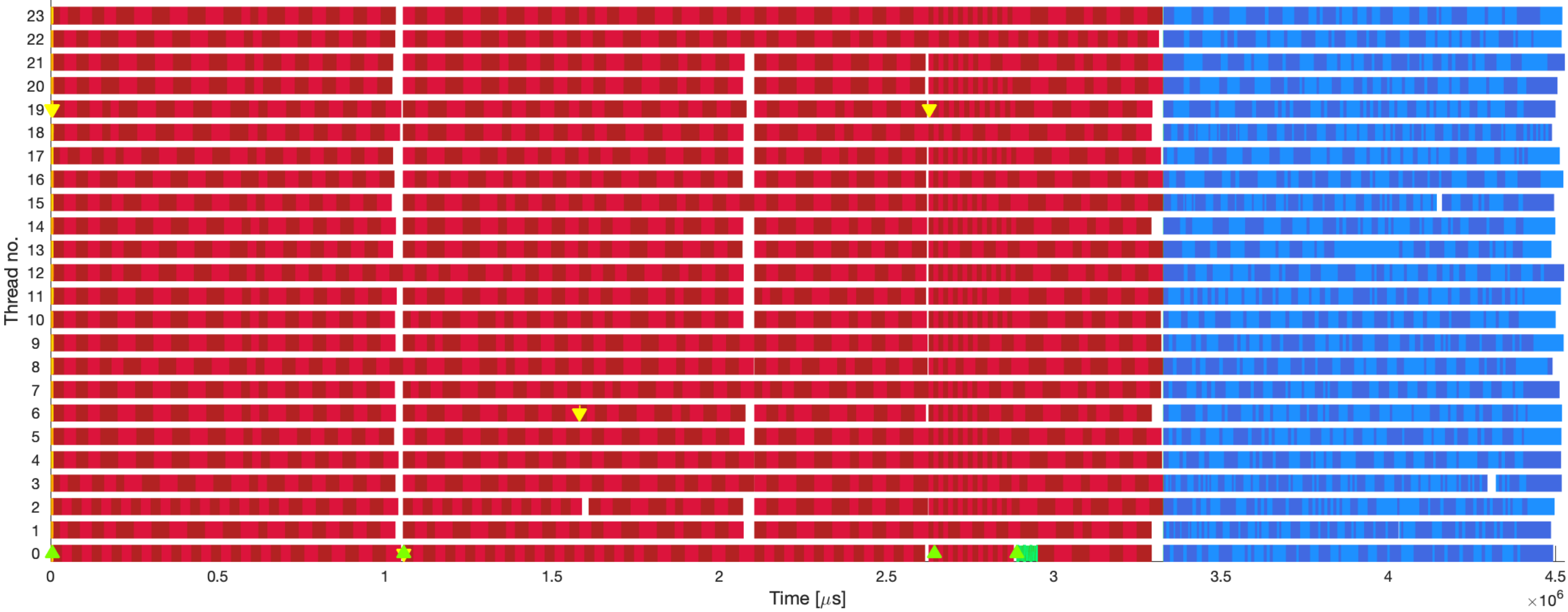}
    \end{center}
    \caption{Execution of the OpenMP tasks on one of Salomon's nodes during a parallel matrix-vector multiplication using 8 nodes. The computation is dominated by the tasks dedicated to the M2L operations (light and dark red rectangles) and nearfield operations (light and dark blue rectangles). The S2M and M2M operations at the beginning of the computation are displayed in orange but are hardly visible (see Figure~\ref{fig:s2m_tasks} for details) and L2L and L2T operations are depicted in green. Finally, moments of MPI communication are marked with yellow (send) and green (receive) triangles. The code was compiled using the Intel Compiler v19.1.1.}
    \label{fig:tasks}
\end{figure}

\begin{figure}[ht]
    \begin{center}
        \includegraphics[width=1\textwidth]{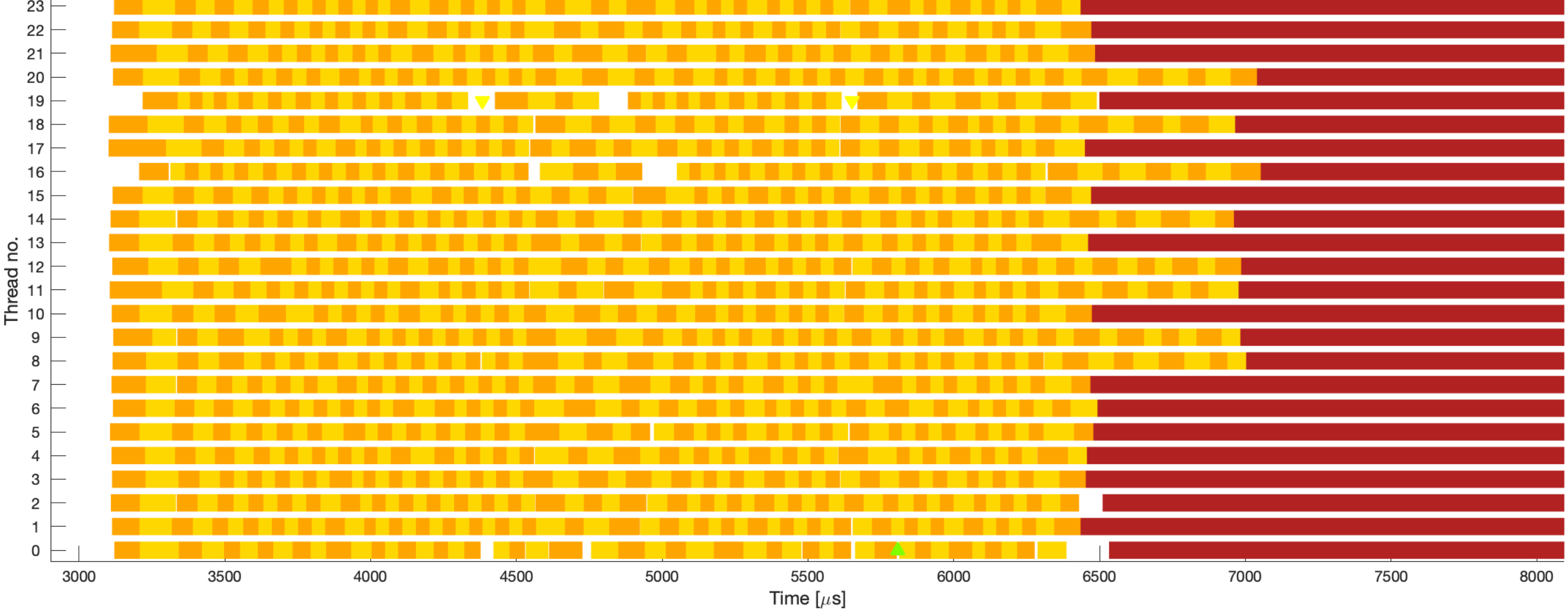}
    \end{center}
    \caption{Details of the S2M and M2M tasks (in orange) from Figure~\ref{fig:tasks}. Moments of MPI communication are marked with yellow (send) and green (receive) triangles.}
    \label{fig:s2m_tasks}
\end{figure}

Next, we solve a Dirichlet problem with 1.5 million space-time boundary elements using eight computational nodes of the Salomon cluster to demonstrate the performance of the task scheduling algorithm presented in Section~\ref{sec:04_parallelization}. We use the same space-time domain as in the previous example discretized into 6144 spatial and 256 temporal elements (equally distributed among 16 time-slices) and the Dirichlet datum given by \eqref{eq:dirdata}. We set $n_{\mathrm{max}}=800, m_t=m_x=6, c_{\mathrm{st}}=0.9$, and $n_{\mathrm{st}}=5$. In Figure~\ref{fig:tasks} we visualize the execution of OpenMP \texttt{tasks} on 24 cores of the node with MPI rank 5. As expected, the computational time is dominated by the \texttt{tasks} dedicated to the M2L computation (denoted by the rectangles in the shades of red) and by the nearfield operations (blue rectangles). L2L and L2T operations are displayed in green and require a negligible amount of the computational time. The S2M and M2M operations at the beginning of computation are depicted in orange but are hardly visible in the graph, therefore we provide a zoom into the first 8000 $\mu \mathrm{s}$ of the computation in Figure~\ref{fig:s2m_tasks}. The moments when data are sent or received via MPI operations are marked with yellow and green triangles, respectively. Note, that the scheduling thread 0 does not only take care of data reception and creating \texttt{tasks} when dependencies are fulfilled, but also participates in the execution of \texttt{tasks}. We observe that the scheduling algorithm (in combination with the Intel OpenMP runtime scheduler) is able to efficiently utilize all available threads. Since the MPI communication is non-blocking, it is hidden by the computation and the scheduling thread can participate in computations and only check for received data whenever it is scheduling new tasks. A very small amount of idle time is still visible in Figure~\ref{fig:tasks}. We could probably overcome these idle times by using a more flexible task scheduler instead of the OpenMP tasks, but since the amount of idle time is negligibly small we kept the latter.

\begin{figure}[b]
    \begin{center}
        \includegraphics[width=1\textwidth]{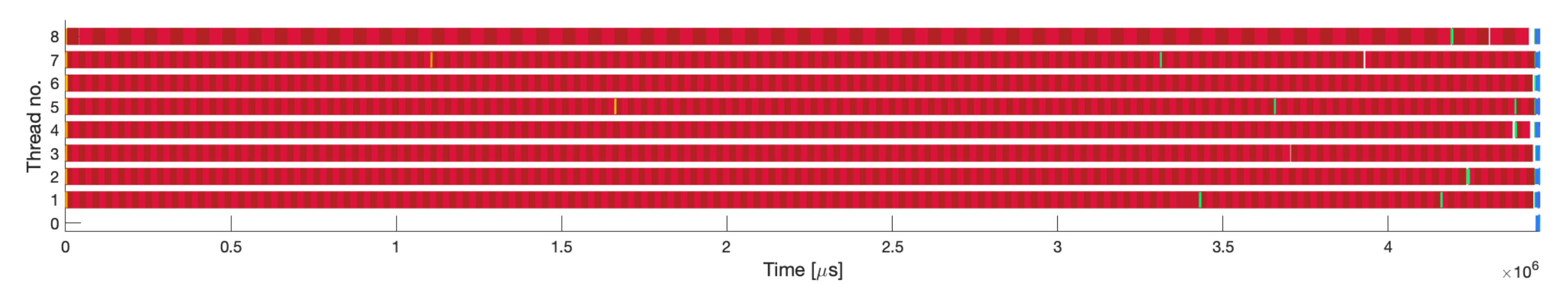}
    \end{center}
    \caption{Exemplary visualization of the OpenMP task execution on one cluster node during a distributed matrix-vector multiplication using GCC v9.3. The \texttt{taskyield} construct shows no effect, and thus the first thread only schedules tasks and does not participate in their execution up until the very end.}
    \label{fig:gcc}
\end{figure}

Compare the results in Figures~\ref{fig:tasks} and \ref{fig:s2m_tasks} with those in Figure~\ref{fig:gcc}. Here we solved a smaller problem on a single node using GCC v9.3 as a compiler. Notice that the first thread does not participate in the \texttt{task} execution up until the very end of the computation since it is busy creating the \texttt{tasks} (and possibly receiving data from other nodes) and the \texttt{taskyield} construct has no effect on it. Of course, this imbalance leads to longer execution times. % A possible remedy could be oversubscribing the number of threads (especially on CPUs with hyperthreading enabled).

 \subsection{Distributed memory performance}

 Distributed memory scalability of the code was tested on up to 256 nodes (6144 cores) of the Salomon cluster. A hybrid MPI-OpenMP parallelization with one MPI process per node and 24 OpenMP threads per process was employed.

 \begin{table}[ht]
    \begin{center}
        \begin{tabular}{rrrrrrrr}
            \toprule
            & \# nodes & 16 & 32 & 64 & 128 & 256 \\
           \midrule
            \multirow{2}{*}{$\matV_{h}$} & time $[s]$ &  769.8 & 385.5 & 194.4 & 97.1 & 50.0 \\
           & efficiency [\%]&  100.0 & 99.8 & 99.0 & 99.1 & 96.2 \\
           \midrule
          \multirow{2}{*}{$\matK_{h}$} & time $[s]$ & 502.4 & 252.5 & 128.6 & 63.0 & 31.9\\
           & efficiency [\%]  & 100.0 & 99.5 & 97.7 & 99.7 & 98.4  \\
           \midrule
          \multirow{2}{*}{iteration} & time $[s]$& 14.7 & 7.3 & 3.7 & 2.1 & 1.5 \\
          & efficiency [\%]  & 100.0 & 101.5 & 99.4 & 89.4 & 62.6 \\
          \bottomrule
        \end{tabular}
    \end{center}
    \caption{Scalability of the code on up to 256 nodes of the Salomon cluster for a problem with 12\,582\,912 space-time surface elements (1024 time-steps, 12\,288 spatial elements).}
    \label{tab:mpi}
\end{table}

We again solve the model problem from Section~\ref{sec:02_bie_bem} with the Dirichlet datum~\eqref{eq:dirdata} and the lateral boundary of the computational domain $(-0.5, 0.5)^3\times (0, 0.25)$ discretized into 12\,288 spatial and 1024 temporal boundary elements (resulting in a total of 12\,582\,912 space-time boundary elements). We equally distribute the temporal elements among 256 time-slices. The variable $n_{\mathrm{max}}$ is set to~800. The orders of the Chebyshev and Lagrange polynomials are set to 12 and~4, respectively, while the parameters $c_{\mathrm{st}}=4.1$ and $n_{\mathrm{tr}}=2$ are used.
Results of the tests are presented in Table~\ref{tab:mpi}. The assembly of the system matrices $\matV_h$ and~$\matK_h$ scales almost optimally up to 256 compute nodes. The same holds for the iteration times but there is a slight drop in efficiency for 256 nodes. This is probably due to the relatively small number of considered time-slices and time-steps. In fact, we decomposed the time interval into 256 time-slices, so when using 256 nodes we assign only one slice to each node and reach the limit of our parallelization scheme for this example. Nonetheless, we achieve a high efficiency for the distributed GMRES solver, reducing the time per iteration from 14.7~s on 16 nodes to 1.5~s on 256 nodes. The iterative solver requires 60 iterations to reach a relative accuracy of $10^{-8}$, thus the total computation time is reduced from 2221~s on 16 nodes to 218~s on 256 nodes.

% The preprocessing phase (i.e., assembly of cluster trees or distributed mesh) is not fully parallelized, therefore it does not scale as well as the remaining parts of the code. However, the preprocessing time remains approximately constant and is currently negligible when compared to the assembly and solve parts of the code.

\begin{figure}[ht]
    \begin{center}
        \includegraphics[width=0.75\textwidth]{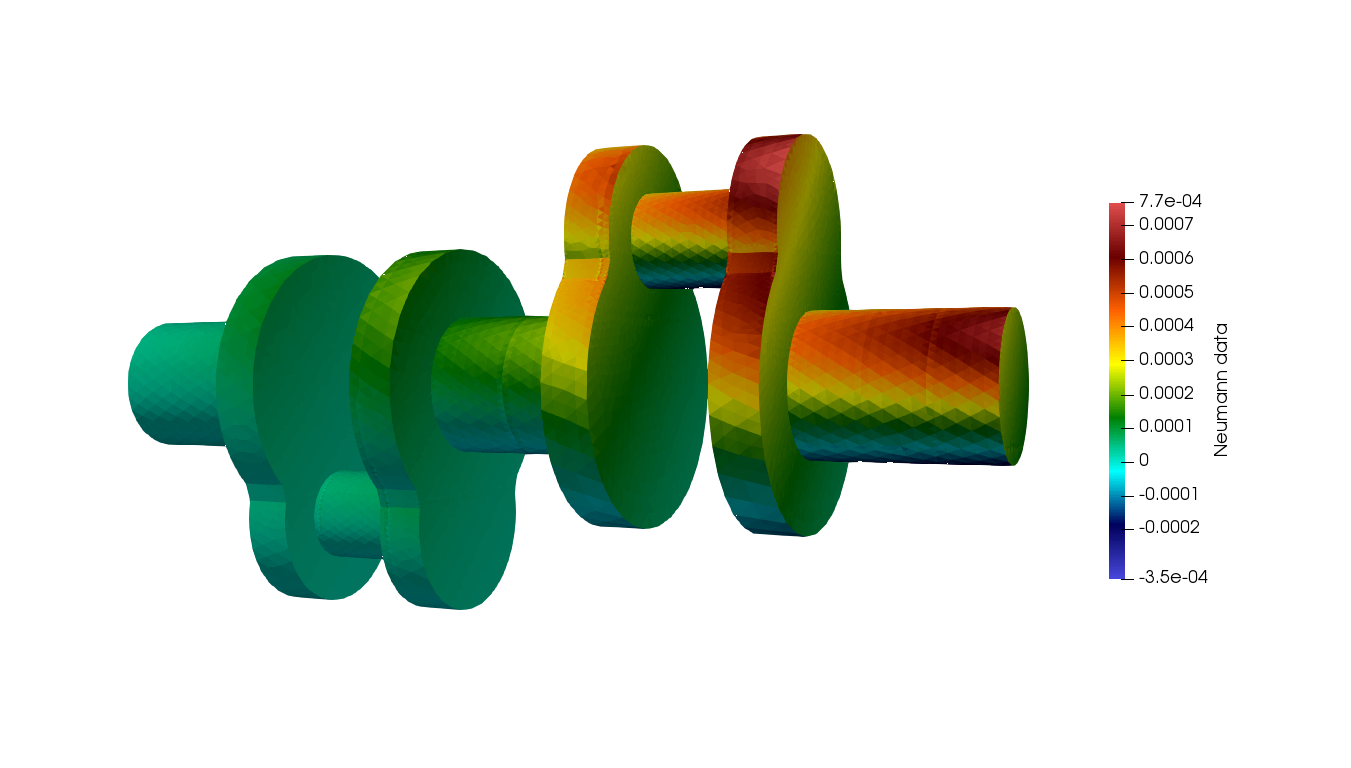}
    \end{center}
    \caption{Computed Neumann datum for the crankshaft discretized by 43\,917\,312 space-time elements at time $t=0.25$.}
    \label{fig:crankshaft}
\end{figure}

Finally, to demonstrate the performance of our code on more realistic examples, we solve the Dirichlet problem for a crankshaft discretized by 42\,888 spatial surface elements and the time interval $(0, 0.25)$ divided into 1024 time-steps (leading to a space-time surface mesh with approximately 44 million boundary elements). We use the same Dirichlet datum as in the previous examples, uniformly distribute the temporal elements among 256 time-slices, and set $n_{\mathrm{max}}=800$, $m_t = 3$, $m_x=12$, $c_{\mathrm{st}}=4.5$, and $n_{\mathrm{tr}}=2$. Using 128 nodes of the Salomon cluster we are able to assemble the system matrices and solve the problem in less than two hours with a relative GMRES accuracy of $10^{-6}$, see Table~\ref{tab:crankshaft}. The solution at the end of the time interval is depicted in Figure~\ref{fig:crankshaft}.

% \begin{table}[ht]
%     \begin{center}
%         \begin{tabular}{rr}
%             \toprule
%             \# compute nodes & 128 \\
%             \midrule
%             $\matV_h$ assembly & 1161.58 s\\
%             $\matK_h$ assembly & 1223.97 s\\
%             \midrule
%             solution & 2927.70 s\\
%             \# iterations & 399\\
%           \bottomrule
%         \end{tabular}
%     \end{center}
%     \caption{Distributed solution of the 'crankshaft problem' with 43\,917\,312 space-time surface elements (1024 time-steps, 42\,888 spatial elements).}
%     \label{tab:crankshaft}
% \end{table}

\begin{table}[ht]
    \begin{center}
        \begin{tabular}{rrrrr}
            \toprule
            \# compute nodes & $\matV_h$ assembly & $\matK_h$ assembly & solution & \# iterations \\
             128 & 1161.58 s & 1223.97 s & 2927.70 s & 399 \\
          \bottomrule
        \end{tabular}
    \end{center}
    \caption{Distributed solution of the ``crankshaft problem'' with 43\,917\,312 space-time surface elements (1024 time-steps, 42\,888 spatial elements).}
    \label{tab:crankshaft}
\end{table}

\section{Conclusion and Outlook}
\label{sec:06_conclusion}
In this paper we developed a parallel space-time FMM for the heat equation. We started from an existing space-time FMM and noticed that its temporal structure can be exploited for parallelization. In fact, the original space-time FMM can be associated with a 1D temporal tree that can be distributed among computing processes and allows to group the actual FMM~operations in time. In our algorithm we used a simple task scheduler to execute these groups of operations in parallel based on individual dependencies and to realize inter-process communication in an asynchronous manner. This allowed us to overcome the strict distinction between FMM phases (forward transformation, multiplication phase, backward transformation) and global synchronization between processes. The resulting algorithm was implemented in the publicly available \cpp{} library BESTHEA~\cite{besthea} as a hybrid MPI-OpenMP code. In several numerical experiments we investigated its efficiency. In particular, we showed close to optimal scalability on a large number of computing nodes.

The parallelism in time of the proposed method is a big advantage over other methods for the solution of boundary value problems of the heat equation like time-stepping schemes. As we have mentioned at the end of Section~\ref{sec:04_parallelization}, one could refine the method by extending the distributed parallelization to an additional spatial dimension. We are optimistic that such an extension is compatible with the presented approach and would allow to increase its scalability even further.

As a next step, we plan to develop a version of the pFMM for adaptive space-time meshes which will enable the development of adaptive BEM for the transient heat equation in 3D.
In addition, we will develop a compression technique for the temporal nearfield blocks suitable for adaptive meshes. This will reduce the storage requirements and make the method more efficient.

\section*{Acknowledgements}
The authors acknowledge the support provided by the Czech Science Foundation under the project 19-29698L, the Austrian Science Fund (FWF) under the project I 4033-N32, and by the Ministry of Education, Youth and Sports of the Czech Republic through the e-INFRA CZ (ID:90140).

\bibliography{references}

\end{document}